\documentclass[english,10pt]{amsart}
\usepackage[english]{babel}
\usepackage{bbm}
\usepackage{hyperref}
\usepackage{verbatim}

\usepackage{amsfonts,amssymb,mathtools}
\usepackage{amsmath}

\usepackage{graphicx,xcolor}
\graphicspath{{Figures/}}
\usepackage{tikz}
\usepackage[style=numeric,maxcitenames=2,maxbibnames=10,doi=false,isbn=false,url=false,natbib=true,giveninits=true,hyperref,bibencoding=utf8,backend=biber]{biblatex}

\newcommand{\N}{\ensuremath{\mathbb{N}}}
\newcommand{\R}{\ensuremath{\mathbb{R}}}

\newcommand{\E}{\ensuremath{\mathbb{E}}}
\renewcommand{\P}{\ensuremath{\mathbb{P}}}

\newcommand{\ind}[1]{\ensuremath{\mathbbm{1}_{\left\{#1\right\}}}}
\newcommand{\diff}{\mathop{}\mathopen{}\mathrm{d}}
\newcommand{\cal}[1]{\ensuremath{\mathcal{#1}}}
\newcommand\croc[1]{\left\langle #1\right\rangle}
\newcommand\steq[1]{\stackrel{\text{\rm #1.}}{=}}

\def\eps{\varepsilon}
\def\cadlag{c\`adl\`ag }
\def\sS{{6S~\mbox{RNA}}}
\def\sSs{{6S~\mbox{RNAs}}}

\newtheorem{proposition}{Proposition}
\newtheorem{definition}[proposition]{Definition}
\newtheorem{lemma}[proposition]{Lemma}
\newtheorem{theorem}[proposition]{Theorem}
\newtheorem{corollary}[proposition]{Corollary}

\title[Stochastic Models of Regulation]{Stochastic Models of Regulation of Transcription in Biological Cells}

\date{\today}
\author[V. Fromion]{Vincent Fromion}
\email{Vincent.Fromion@inrae.fr}
\address[V.~Fromion, J. Zaherddine]{INRAE, MaIAGE, Université Paris-Saclay, Domaine de Vilvert, 78350 Jouy-en-Josas, France}
\author[Ph. Robert]{Philippe Robert}
\email{Philippe.Robert@inria.fr}
\urladdr{http://www-rocq.inria.fr/who/Philippe.Robert}
\author[J. Zaherddine]{Jana Zaherddine}
\email{Jana.Zaherddine@inria.fr}
\address[Ph.~Robert, J. Zaherddine]{INRIA Paris, 2 rue Simone Iff, 75589 Paris Cedex 12, France}

\begin{document}

\begin{abstract}
In this paper we study an important {\em global } regulation mechanism of transcription of biological cells using specific macro-molecules, {\sSs}. The functional property of  {\sSs} is of blocking  the transcription of RNAs when the environment of the cell is not favorable. We investigate the efficiency of this mechanism with a scaling analysis of a stochastic model. 
The evolution equations of  our  model are driven by the law of mass action and the total number of polymerases is used as a scaling parameter. Two regimes are analyzed: exponential phase when the environment of the cell is favorable to its growth, and the stationary phase when resources are scarce. In both regimes, by defining properly occupation measures of the model,  we prove an averaging principle for the associated  multi-dimensional Markov process on a convenient timescale, as well as convergence results for ``fast'' variables of the system. An analytical expression of the asymptotic fraction of sequestrated polymerases in stationary phase is in particular obtained. The consequences of these results are discussed. 
\end{abstract}

\maketitle

 \vspace{-5mm}

\bigskip

\hrule

\vspace{-3mm}

\tableofcontents

\vspace{-1cm}

\hrule

\bigskip
\section{Introduction}
The central dogma of molecular biology asserts for biological cells that genetic information flows mainly in one direction, from DNA to RNAs, and to proteins.
For the two most studied bacteria {\it Escherichia coli} and {\it Bacillus subtilis}, production of  proteins is a central process which can be described as a process in two main steps. In the first step, macro-molecules {\em polymerases} produce RNAs with genes of DNA. This is the {\em transcription} step. The second step, {\em translation}, is the production of proteins from mRNAs, {\em messenger RNAs}, with macro-molecules {\em ribosomes}. See~\citet{Watson}.

In bacterial cells,  protein  production uses essentially most of cell resources: a large number of its macro-molecules such as polymerases and ribosomes,  biological bricks of proteins, i.e.,  amino acids, and  the energy necessary to build proteins during the translation step, such as GTP.

In this paper we study an important regulation mechanism of transcription using specific RNA macro-molecules, {\sSs},  common to a large number of bacteria. See~\citet{Wassarman2018} for example. The functional property of this RNA is of blocking/sequestering  free polymerases from producing RNAs. The general context of this regulation is related to complex mechanisms of the cell to finely tune the production of a large set of RNAs. Let us first recall the three main categories of RNAs:
\begin{enumerate}
\item RNAs used for the building of ribosomes, i.e., rRNAs, {\em ribosomal RNAs}.  A ribosome is a complex assembly of around 50 proteins and, also, of several rRNAs. An rRNA is a long chain of several thousands of nucleotides, it is in particular a costly  macro-molecule to produce. Reducing or speeding-up the production of ribosomes, in particular of rRNAs, has therefore a critical  impact on resource management of the cell;
\item mRNAs, {\em messenger RNAs}, used by the translation step to produce a protein from mRNAs coding sequences;
\item  A large set of RNAs that do not belong to the two previous categories, such as  {\em transfer RNAs}, tRNAs, or  {\em Bacterial small RNAs}, sRNAs, often associated to  regulation mechanisms. This class  includes  {\sSs}.
\end{enumerate}
When the concentrations of different resources in the medium  are high enough  for some time, the bacterium has the ability to use them efficiently, via its complex regulatory system, to  reach a stable  exponential growth regime with a fixed  growth rate. The growth of a bacterial population in a given medium leads therefore to an active consumption of resources necessary to produce new cells.

When resources are scarce, each bacterium of the population can adapt,  to either exploit differently the available resources, or to do without some of them, as for example when some amino-acids are missing. For {\it E. coli} or {\it B. subtilis}, these bacteria  use in priority  resources maximizing their growth rate. In the context of this adaptation, and for reasons related to the decay of resources,  each bacterial cell  has to decrease its growth rate, and finally to ultimately stop its growth.

 The regulatory network involved in the management of the growth rate to adapt to the environment is complex. The important point in this domain is that the bacterium has to modify the concentration of most of the agents  in charge of it:  number of ribosomes, concentrations of proteins in the metabolic network,   transporters, \ldots In a first, simplified, description,  the decay of a specific resource in the environment leads  to a move to a state of the cell where  concentrations of several components have been adapted. To study  the transition between growth phases, we have chosen to focus on the action of a small RNA, {\sS},  which  plays an important, even essential, role in this domain. Note that, even if this mechanism is central, this description of the transition between  growth regimes is nevertheless a simplification in our approach, since the bacterial cell has different ways to modify the steady-state level of its components.

 In this article, we investigate a simplified scenario where transitions occur between two phases:  an exponential growth phase and  the stationary phase, where the growth rate is equal to $0$. The first interest of this scenario lies in the sharp transition of the polymerase management by the cell, via the strong effect of the {\em stringent response} on the production of rRNAs: the transcription of rRNAs is completely stopped. This  is where the action of {\sSs} is crucial. See~\citet{Gottesman2006}. Its second interest is experimental since its is possible in practice to create this transition by the addition of a convenient product in the medium of  cell populations to induce a stringent response.
  Our general goal is to investigate if, with this simplified framework,  the regulatory system organized around {\sS}  has the desired qualitative properties to ensure a convenient transition between these regimes.  In this paper,  we analyze the efficiency of the regulation by {\sSs}  with stochastic models. We investigate in particular the time evolution of the activity of polymerases in the cell under different regimes.
 
\subsection{A Simple Description as a Particle System}
In order to explain the basic principle of the  regulation mechanism investigated in this paper, we describe a simplified version in terms of a particle system. Section~\ref{BioSec} describes in more depth and detail the biological context of this class of models. 

We consider two types of particles $P$ and $6S$. There is a fixed number of particles of type $P$ and there are random arrivals of particles of type $6S$. A particle of type $P$ can be in three states: busy, idle, or paired with a $6S$ particle. Similarly, a $6S$ particle is either idle or paired. The possible events are:
\begin{itemize}
\item an idle, resp. busy,  $P$ particle becomes busy, resp. idle;
\item a couple of an idle $P$ and an idle $6S$ is paired;
\item a pair $P{-}6S$ is broken giving two idle $P$ and $6S$;
\item an idle $6S$ arrives/dies.
\end{itemize}
Note that only an idle $6S$ can die. A statistical assumption is that each couple of free $P$ particle and free $6S$ particle is paired at some fixed rate and each free $6S$ dies at a fixed rate too.

We present a heuristic description of the phenomena we are interested in: 
\begin{enumerate}
\item\label{itema} If the parameters of the $P$ particles are such that, on average,  most of particles of type $P$ are busy. Therefore, few of them are idle, the arriving $6S$ particles will very likely die before they can be paired with a $P$ particle. In this case there will be few $6S$ particles in the system.
\item\label{itemb} Otherwise, if, on average,  a significant fraction of particles of type $P$ are idle, the arriving $6S$ particles will very likely pair with one of them. In particular, as long as there are many idle $P$ particles, $6S$ will be quickly paired so that few of them will die. In this manner, the dynamic arrivals of $6S$  progressively decrease the number of  idle $P$ particles.
\end{enumerate}
A pair $P{-}6S$ is seen as a sequestration of a $P$ particle, the purpose of $6S$ particles is of storing ``useless'' $P$ particles.  The case~\ref{itema}) corresponds to the case when most of $P$ particles are efficient so that no regulation is required. This corresponds to the {\em exponential growth phase} of our biological process. The case~\ref{itemb}) is when there is a need of sequestration of $P$ particles, this is the {\em stationary growth phase} of our model.

The nice feature of this mechanism is its adaptive property due to the dynamic arrivals of $6S$: if they are useless, they disappear after some time. Otherwise, as it will be seen, their number builds up until some threshold of sequestration is reached. 

The main goal of the paper is of understanding under which conditions on the parameters the cases~\ref{itema}) or~\ref{itemb}) may occur.  To assess the efficiency of the regulation mechanism in the case~\ref{itema}), we study the time evolution of the number of $6S$  particles. In the case~\ref{itemb}),  we  investigate the number of sequestered $P$ particles to determine the maximal sequestration rate of the regulation.

The model investigated in the paper is in fact a little more complicated in the sense that $P$ particles can be ``busy'' in two ways: either it remains busy during a random amount of time before being idle again.  The other busy state is that it joins a  queue where only the particle at the head of the queue becomes idle again after a random amount of time. In our model, this is related to mRNAs and rRNAs production. The next section gives a detailed description of these aspects. 

\subsection{Biological Background}\label{BioSec}
\subsubsection*{Transcription}
In a bacterial cell, a polymerase may be associated to several specific proteins, called $\sigma$-factor to form a {\em holoenzyme} E$\sigma$. In our case we focus on the ``housekeeping'' $\sigma$-factor $\sigma^{70}$.  This holoenzyme binds to a large set of gene promoters to initialize the transcription. This is the {\em initiation phase}. If this step is successful, the protein $\sigma^{70}$ is detached and the polymerase completes the elongation of the corresponding RNA.

This is a simplified description of course. The precise description of the  mechanisms  are dependent on the bacterium,  it is nevertheless sufficiently accurate from our modeling perspective. Throughout the paper we do the slight abuse of using the term polymerase instead of the more biologically correct term holoenzyme.  Another important aspect is that the initiation phase may fail due to random fluctuations within the cell, or to a low level of nucleotides needed for the initiation of transcription, i.e. ATP, GTP, UTP, CTP,~\ldots When this happens the transcription is aborted.   The level of GTP, for example,  has an impact on the modulation of initiation of transcription with respect to the growth rate for bacterium B. subtilis,  and, similarly, the level of ppGPpp for bacterium \textit{Escherichia coli}. See~\citet{Geissen} in the case of an rRNA. 
\subsubsection*{Regulation by small RNAs}
A subset of RNAs whose sizes in nucleotides is less than $100$, small RNAs or sRNAs  has been shown to play an important role to regulate gene expression. The first such sRNAs were identified in the late 1960's. See~\citet{Britten} and~\citet{Zamore}.  They were shown to turn in or turn off specific genes under convenient conditions.

Among them the sRNA {\sS} was first discovered  because of its abundance in E. coli  in some circumstances. See~\citet{Hindley}. This has been one of the first sRNAs to be sequenced.  Nevertheless, it took three decades to understand its role in the regulation of transcription.

Experimental studies have shown that {\sS} acts in fact as a {\em global regulator} of  transcription and not only for the regulation of a reduced subset of genes as most of small RNAs. A  {\sS} has a three-dimensional structure similar to a DNA promoter, so that the holoenzyme” E$\sigma^{70}$ may be bound to it and is, in some way, sequestered by it.  See \citet{Cavanagh2014} and~\citet{Nitzan2014}.   It has been shown that  during {\em stationary phases}, when the growth rate is null, the {\sSs} accumulate to a high level, with more than $10000$ copies. During an {\em exponential phase}, when the growth is steady, the average duration time of cell division is around 40min for E. coli, there are less than $1000$ copies. See~\citet{Wassarman2018} and~\citet{Steuten2014}.

The fluctuations of the number of copies of {\sSs} is therefore an important indicator of the growth rate of the cell. An important question is to assess the efficiency of the regulation mechanism operated by the {\sSs}. The impact of several parameters of the cell are investigated: The total number of polymerases, the production rate of {\sS} and their degradation rate, initiation rates of polymerases for rRNAs and mRNAs and the sequestration rate, i.e. the binding rate of a couple {\sS} and polymerase.
    
\subsection{Mathematical Models}
Regulation of gene expression has been analyzed with mathematical models for some time now. See~\citet{Mackey} and also Chapter~6 of~\citet{Bressloff}, and the references therein. The lac operon model is  one of the most popular mathematical models in this domain, for its bistability properties in particular.  See also~\citet{DFR}.

Specific stochastic models of regulation by RNAs are more scarce.  The regulation of mRNAs by sRNAs in a stochastic framework has been the subject of several studies recently. In~\citet{Kumar}, \citet{Mehta}, and~\citet{Platini}, the authors study regulation mechanisms of mRNAs by sRNAs with a two-dimensional Markov chain for the time evolution of the number of sRNAs and mRNAs.  Some limiting regimes of  the corresponding Fokker-Planck evolution equations are investigated and discussed. The difficulty is the quadratic dependence on the number of mRNAs and sRNAs. See also~\citet{Baker2012} and~\citet{Mitarai}. These studies can be seen as extensions of the early works on stochastic models of gene expression, see~\citet{Berg}, \citet{Elowitz} and~\citet{Rigney}.

\subsection{The Main  Results}
In this paper, we will study the efficiency of the sequestration of polymerases by {\sSs}. Recall that this is in fact the holoenzyme which is sequestered. We investigate the  behavior of several variables associated to the regulation of the transcription phase: Number of free/sequestered polymerases and number of polymerases in the elongation phase of mRNAs and rRNAs.

\subsection*{Technical Challenges}
We assume that there are $N$ polymerases with $N$ large. We derive functional limiting results, with respect to this scaling parameter, of the time evolution of several stochastic processes. An important feature of our model is that the main Markov process exhibits a quadratic dependence of the state of the process, due to the use of the law of mass action for the dynamic of our model.   One of  the main technical difficulties is in the proof of Theorem~\ref{OLLN} of an averaging principle. Several preliminary results have to be established as well as a convenient definition of occupation measures. This is  due to the (very) fast underlying timescale, $t{\mapsto}N^2 t$,   used.  Formally, the diffusion component is of the order of $N$ but should vanish for this first order result.  For this reason, in a first step, the ``slow'' processes are included in the definition of occupation measures and not only the ``fast'' processes as it is done in general. In our proofs we use several coupling arguments,  estimates of hitting times of rare events, stochastic calculus for stochastic differential equations driven by Poisson processes, and the framework of averaging principles.

\subsection*{A Chemical Reaction Network Description}
 For simplicity this number is assumed to be constant. There is also a production of {\sSs} which we will distinguish from the production of other RNAs. From the point of view of our model, polymerases can be in several states
\begin{itemize}
\item {\em Free.} The polymerase may bind to a gene of an mRNA, or of an rRNA, or be sequestered by a {\sS}, $(F_N(t))$ denotes the process of the number of  free polymerases. 
\item {\em Transcription of an mRNA. }
A  chain of nucleotides is produced, $(M_N(t))$ is the number of such polymerases. There is a large number of types of mRNAs.
\item {\em Transcription of an rRNA.}
A long chain of nucleotides is produced. As it will be seen, it is described by a  process $((U_j^N(t),R_j^N(t)), 1{\le}j{\le}J)$. The  number $J$ of types of rRNAs is usually small, less than ten. We denote by $\|R(t)\|$ the total number of polymerases in this situation.  
\item {\em Sequestered by a {\sS}. } The associated process is $(S_N(t))$.
\end{itemize}
Similarly a {\sS} can be either free or paired with a polymerases, $(Z_N(t))$ denotes the process of the number of free {\sSs}.
See Section~\ref{ModSec} for more details.

With these notations,  the assumption on the conservation of mass for the polymerases gives the relation
\[
F_N(t){+}M_N(t){+}S_N(t){+}\|R_N\|(t){=}N,\qquad \forall t{\ge}0.
\]
The dynamic of this stochastic system is governed by the analogue of the law of the mass action in this context. See~\citet{AndersonKurtz}. The rate of creation of sequestered polymerases is in particular quadratic with respect to the state, it is proportional to $F_NZ_N$. This is one of the important features of this stochastic model.

\subsection*{Two Limiting Regimes}
Our mathematical results can be described  as follows. See the formal statements in Section~\ref{ExpSec} and~\ref{StatSec}. In Definition~\ref{DefPhase}, we introduce two sets of conditions on the parameters of our model, which define the exponential regime and the stationary regime,  our cases a) and b) above. We do not detail them here. Assuming  that the maximum number of polymerases simultaneously in transcription of rRNAs, resp. mRNAs, is of the order of $c_rN$, resp. $c_mN$,   under some scaling conditions and appropriate initial conditions, we have: 
\begin{enumerate}
\item[1)] {\em Exponential Phase.}\\
 For the convergence in distribution
  \[
  \lim_{N\to+\infty} \left(\frac{\|R_N\|(t)}{N},\frac{M_N(t)}{N}\right)=(c_r,1{-}c_r),
  \]
  and, for any  $t_0{>}0$, the random variable $(F_N(t_0))$ converges in distribution to a Poisson distribution and
  the sequence of process $(S_N(t),Z_N(t))$ is converging in distribution to a positive recurrent  Markov process.
  See Theorem~\ref{TheoExpFree}.

In this case, the polymerases are mostly doing transcription of rRNAs or mRNAs, few of them are free or sequestered by a {\sS}. 
\item[2)] {\em Stationary Phase.}\\
   For the convergence in distribution
  \[
  \lim_{N\to+\infty} \left(\frac{M_N(Nt)}{N}, \frac{F_N(Nt)}{N},\frac{S_N(Nt)}{N}\right)=(c_m,\overline{f}(t),1{-}c_m{-}\overline{f}(t)),
  \]
  where $(\overline{f}(t))$ is the solution of an ODE, such that
  \[
  \lim_{t\to+\infty}\overline{f}(t)= \overline{f}(\infty)>0.
  \]
The process $(\|R_N(t)\|)$ is stochastically upper-bounded by a positive recurrent Markov process. See Theorem~\ref{theostat}.

In the stationary phase  there are few polymerases doing transcription of rRNAs. A fraction of them remains free, asymptotically $\overline{f}(\infty)$, i.e. the sequestration process does not control all ``useless'' polymerases. This is in fact a non-trivial consequence of  the dynamic creation and destruction of {\sSs}, even if an {\sS} paired with a polymerase cannot be degraded. The fact that sequestration phenomenon of a fraction of the $N$ polymerases occurs on the time scale $t{\mapsto}Nt$ is intuitive given that the rate of creation of {\sSs} is constant. 
\end{enumerate}
In all cases the process $(Z_N(t))$ is stochastically upper-bounded by a positive recurrent Markov processes.

\subsection{Outline of the Paper}
Section~\ref{SecMod} introduces in detail the complete model of transcription and also an important model, the auxiliary process. The exponential/stationary phases corresponds to super/sub critical condition for this auxiliary process. They are investigated respectively in Section~\ref{SecSub} and~\ref{SecSuper}. The last  sections~\ref{ExpSec} and~\ref{StatSec} are devoted to the exponential/stationary regimes of our model. 
\section{Stochastic Model}\label{SecMod}
In this section we introduce the state space description of the regulation of transcription. We first describe our main assumptions in the design of the stochastic model.
\subsection{Modeling Assumptions}\label{ModSec}
The chemical species involved in the regulation  process are the genes of different types of mRNAs and rRNAs and of $\sS$, and polymerases. The products are  different types of mRNAs and of rRNAs and also {\sSs}. 
\begin{itemize}
\item {\sc Transcription of }rRNAs.\\
 There are $J$ types of rRNAs and there is a promoter (binding site for polymerases) for each of them. 
 The transcription of an rRNA of type $j$, $1{\le}j{\le}J$,  is in two steps. The promoters of rRNAs are assumed to have a high affinity during the growth phase: If one of these promoters is empty and if there is at least one free polymerase, then the promoter is occupied right away by a polymerase.
 
Once a polymerase is bound to the promoter of the rRNA of type $1{\le}j{\le}J$, it starts elongation at rate $\alpha_{r,j}$  if there are strictly less than $C_{r,j}^N$ polymerases in the elongation phase of this rRNA.  At a given moment there cannot be more than $C_{r,j}^N$ polymerases in elongation of an rRNA of type $j$. 

For each polymerase in elongation,  nucleotides are collected at rate $\beta_{r,j}$. The simplification of the model is that the polymerases in elongation are moving closely on the gene so that  the duration of time to get the last nucleotide for the oldest polymerase in elongation is enough to describe the time evolution of the number of polymerases producing  rRNA of type $j$. 
The  polymerases associated to an rRNA of type $j{\in}\{1,\ldots,J\}$ can then be  represented as a couple $(u_j,R_j)$, where $u_j{\in}\{0,1\}$ indicates if a polymerase is on the promoter or not, and $R_j{\in}\N$ is the number of polymerases in elongation: If $R_j{\ge}1$, an rRNA of type $j$ is therefore created at rate $\beta_{r,j}$.

The assumption is reasonable in the exponential phase, since in this case the number of polymerases producing rRNA of type $j$ is maximal, of the order of $C_{r,j}^N$.  See Section~\ref{ExpSec}. The rRNA part of the system is therefore saturated. In the stationary phase, this assumption has little impact since, as we shall see, the total number of polymerases in the elongation phase of rRNAs is small with high probability and, therefore, negligible for our scaling analysis.

\item[]
\item {\sc Transcription of} mRNAs.\\
It is assumed that there are $C_m^N$ different types of mRNAs and that at a given time, for any $1{\le}i{\le}C_m^N$ there is at most one polymerase in the elongation phase of an mRNA of type $i$. When the promoter of an mRNA of type $i$ is free, a free polymerase may bind to this promoter at a rate $\alpha_m$.  If the promoter of an mRNA of type $i$ is occupied, an mRNA is released at rate $\beta_m$ and the corresponding polymerase leaves the promoter at that instant. The production of mRNAs have simplified in the sense that the initiation phase and the elongation phase are merged into one step. The results obtained in this paper could be obtained without too much difficulty for a model distinguishing them, but at the expense of a more complicated state variable.

The main difference in our model  between the rRNAs and the mRNAs is on the number of polymerases in elongation of the corresponding gene. At a given moment, under favorable growth conditions, there will be many polymerases in the elongation phase of an rRNA,  due in particular to the high initiation rate of these genes. For the  mRNAs, the number of polymerases in elongation phase of a given mRNA type should be small in general. Indeed, there are in each cell few copies of each messenger (from 1 to 100).  Furthermore, the rate of production of each messenger is such that its small number remains  on average constant during  growth or stationary phases and despite the regular degradation (average of  2 minutes half-life in  high-growth rate phase) of each of them.  See Section~\ref{NumbSec}. In our model we have set the  maximal number of polymerases in elongation phase of a given mRNA type  to $1$ for simplicity, but it is not difficult to adapt our results with a maximum number $D$.  Similarly, the initiation rates and production rate, $\alpha_m$ and $\beta_m$ are taken equal for all species of mRNA, also for the sake of simplicity. We have simplified the description of the production of mRNAs to focus  mainly on the sequestration mechanisms that regulate the transcription. From our point of view, the production of mRNAs  holds/stores a subset of polymerases and releases each of them after some random amount of time. It should be noted that this is in fact the usual mathematical setting  to investigate the fluctuations of the production of mRNAs and proteins.  See~\citet{Berg} and~\citet{Rigney}, see also~\citet{Paulsson} for a review of these models.

\item[]
\item {\sc Creation/Degradation of} {\sSs}.\\
  The creation of {\sSs} involves, of course, polymerase. As in the case of mRNAs, it is assumed that there is at most one polymerase in elongation phase of this sRNA. A {\sS} is created at rate $\beta_6{>}0$. A {\sS} is  free when it is not bound to a polymerase. A given free {\sS} is degraded at rate $\delta_6{\ge}0$. Only a free {\sS} can be degraded. 
\item[]
\item {\sc Sequestration/de-Sequestration of}  Polymerases.\\
A polymerase is free when it is not bound to a gene or to a {\sS}. In our study the total number of polymerases is assumed to be constant equal to $N$. A free polymerase is bound to a free {\sS} at rate $\lambda$.  A complex polymerase-{\sS} breaks into a free polymerase and a free {\sS} at rate $\eta$.
\end{itemize}

\subsection{The Markov Process and its $Q$-Matrix}\label{QMatrix}
The vector $(\alpha_{r,j})$ introduced is the vector of {\em initiation} rates of transcription of the different types of rRNAs. The difference between a slow growth (stationary phase) and a steady growth (exponential phase) will be expressed in terms of the comparison, coordinate by coordinate, of the vectors $(\alpha_{r,j})$  and $(\beta_{r,j})$.  We now give a Markovian description of our system. Convenient limiting results will be obtained for the associated Markov process in both phases. 

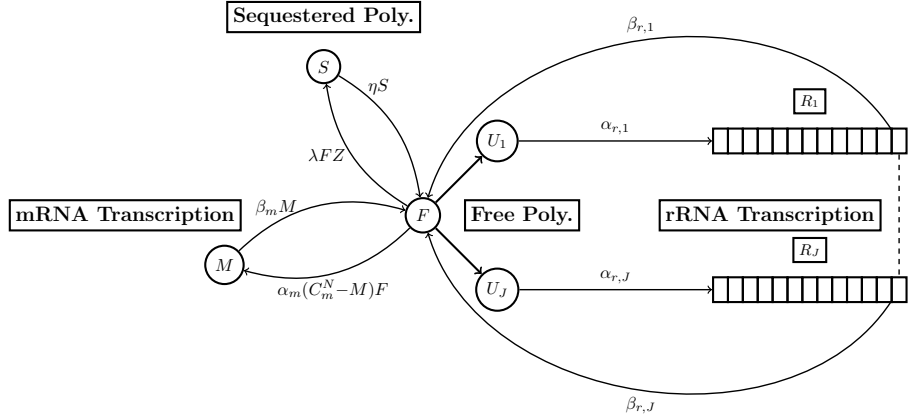
\begin{figure}[ht]
\resizebox{12.0cm}{!}{%
\begin{tikzpicture}[->,node distance=8mm]

  \node[black,very thick,circle,draw](F) at (4,0){$F$};
  \node [black,very thick,rectangle,draw] at (6,0) {\Large\bf Free Poly.};
  \node[black,very thick,circle,draw](U_1) at (5.5,1.5){$U_1$};
  \node[black,very thick,circle,draw](U_J) at (5.5,-1.5){$U_J$};
  \node[black,very thick,rectangle,minimum width = 0.3cm, minimum height = 5mm,draw](R_11) at (10,1.5){};
  \node[black,very thick,rectangle,minimum width = 0.3cm, minimum height = 5mm,draw]() at (10.3,1.5){};
  \node[black,very thick,rectangle,minimum width = 0.3cm, minimum height = 5mm,draw]() at (10.6,1.5){};
  \node[black,very thick,rectangle,minimum width = 0.3cm, minimum height = 5mm,draw]() at (10.9,1.5){};
  \node[black,very thick,rectangle,minimum width = 0.3cm, minimum height = 5mm,draw]() at (11.2,1.5){};
  \node[black,very thick,rectangle,minimum width = 0.3cm, minimum height = 5mm,draw]() at (11.5,1.5){};
  \node[black,very thick,rectangle,minimum width = 0.3cm, minimum height = 5mm,draw](R_12) at (11.8,1.5){};
  \node[black,very thick,rectangle,minimum width = 0.3cm, minimum height = 5mm,draw]() at (12.1,1.5){};
  \node[black,very thick,rectangle,minimum width = 0.3cm, minimum height = 5mm,draw]() at (12.4,1.5){};
  \node[black,very thick,rectangle,minimum width = 0.3cm, minimum height = 5mm,draw]() at (12.7,1.5){};
  \node[black,very thick,rectangle,minimum width = 0.3cm, minimum height = 5mm,draw]() at (13,1.5){};
  \node[black,very thick,rectangle,minimum width = 0.3cm, minimum height = 5mm,draw]() at (13.3,1.5){};
  \node[black,very thick,rectangle,minimum width = 0.3cm, minimum height = 5mm,draw](R_13) at (13.6,1.5){};
  \node (R_1) [above of=R_12,black,very thick,rectangle,draw]{{\small\bf $R_1$}};

  \node[black,very thick,rectangle,minimum width = 0.3cm, minimum height = 5mm,draw](R_J1) at (10,-1.5){};
  \node[black,very thick,rectangle,minimum width = 0.3cm, minimum height = 5mm,draw]() at (10.3,-1.5){};
  \node[black,very thick,rectangle,minimum width = 0.3cm, minimum height = 5mm,draw]() at (10.6,-1.5){};
  \node[black,very thick,rectangle,minimum width = 0.3cm, minimum height = 5mm,draw]() at (10.9,-1.5){};
  \node[black,very thick,rectangle,minimum width = 0.3cm, minimum height = 5mm,draw]() at (11.2,-1.5){};
  \node[black,very thick,rectangle,minimum width = 0.3cm, minimum height = 5mm,draw]() at (11.5,-1.5){};
  \node[black,very thick,rectangle,minimum width = 0.3cm, minimum height = 5mm,draw](R_J2) at (11.8,-1.5){};
  \node[black,very thick,rectangle,minimum width = 0.3cm, minimum height = 5mm,draw]() at (12.1,-1.5){};
  \node[black,very thick,rectangle,minimum width = 0.3cm, minimum height = 5mm,draw]() at (12.4,-1.5){};
  \node[black,very thick,rectangle,minimum width = 0.3cm, minimum height = 5mm,draw]() at (12.7,-1.5){};
  \node[black,very thick,rectangle,minimum width = 0.3cm, minimum height = 5mm,draw]() at (13,-1.5){};
  \node[black,very thick,rectangle,minimum width = 0.3cm, minimum height = 5mm,draw]() at (13.3,-1.5){};
  \node[black,very thick,rectangle,minimum width = 0.3cm, minimum height = 5mm,draw](R_J3) at (13.6,-1.5){};
  \node (R_J) [above of=R_J2,black,very thick,rectangle,draw]{{\small\bf $R_J$}};
  \node () [black,very thick,rectangle,draw]  at (11,0) {{\Large\bf rRNA Transcription}};
  \path (F) edge [black,very thick] (U_1);
  \path (F) edge [black,very thick] (U_J);
  \path (U_1) edge [black,thick,left,midway,above] node {$\alpha_{r,1} $} (R_11);
  \path (U_J) edge [black,thick,left,midway,above] node {$\alpha_{r,J} $} (R_J1);
  \path (R_13) edge[black,thick,dashed,right,-] node{} (R_J3);
  \path (R_13) edge [black,thick,bend right=65,midway,above]  node {$\beta_{r,1}$} (F);
  \path (R_J3) edge  [black,thick,bend left=65,midway,below]  node {$\beta_{r,J} $} (F);

  \node[black,very thick,circle,draw](M) at (0,-1){$M$};
  \path (F) edge [black,thick,bend left,midway,below] node {$\alpha_m(C_{m}^N{-}M)F $} (M);
  \path (M) edge [black,thick,bend left,near start,above] node {$\beta_m M $} (F);
  \node () [black,very thick,rectangle,draw]at (-2,0){{\Large\bf mRNA Transcription}};

  \node[black,very thick,circle,draw](S) at (2,3){$S$};
  \path (F) edge [black,thick,bend left=25,midway,left] node {$\lambda F Z $} (S);
  \path (S) edge [black,thick,,bend left=25,very near start,right] node {$\ \eta S$} (F);
  \node () [black,very thick,rectangle,draw] at (2,4){{\Large\bf Sequestered Poly.}};

\end{tikzpicture}%
}%
\caption{Polymerases: Transcription of mRNAs/rRNAs and Sequestration}
\end{figure}

\subsection*{State Space}
All transitions described in the last section occurs after a random amount of time with an exponential distribution. With this assumption, there is a natural Markov process to investigate the regulation of transcription. The state space is given by
\begin{multline*}
  {\cal S}_N\steq{def}\left\{\rule{0mm}{6mm}x{=}(f,s,z, (u_j,r_j)){\in}\N^3{\times}\prod_{j=1}^J\left(\rule{0mm}{4mm} \{0,1\}{\times}\{0,\ldots,C_{r,j}^N\}\right):\right.\\
  \left. f{+}s{+}\sum_{j=1}^J(u_j{+}r_j)\le N\text{ and if } f{>}0, \text{ then }u_j{=}1, \forall 1{\le}j{\le}J \right\},
\end{multline*}
if the state of the system is $x{=}(f,s,z, (u,r)){\in}{\cal S}_N$, then
\begin{itemize}
    \item $f$ is the number of free polymerases;
    \item $s$, the number of sequestered polymerases;
    \item $z$, the number of free {\sSs};
      \item[] $(u,r){=}((u_j,r_j),1{\le}j{\le}J)$,
    \item $u_{j}{\in}\{0,1\}$ to indicate if a polymerase is bound to the promoter of the  rRNA  of type $j$ or not;
    \item $0{\le}r_{j}{\le}C_{r,j}^N$,  number of polymerases in elongation phase of an rRNA of type~$j$.
    \item In state $x$, the number of polymerases in elongation phase of an mRNA is given by
\begin{equation}\label{MmRNA}
      \Psi(x) \steq{def} N{-}f{-}s{-}\sum_{j=1}^J(u_j{+}r_j).
\end{equation}
\end{itemize}

The associated Markov process is denoted by 
\[
(X_N(t))\steq{def} \left(F_N(t),S_N(t),Z_N(t),(U_N(t),R_N(t))\right),
\]
with $(U_N(t),R_N(t)){=}((U^N_{j} (t),R^N_{j}(t)),1{\le}j{\le}J)$.
The number of polymerases at time $t$ in elongation phase of an mRNA is defined by $M_N(t){=}\Psi(X_N(t))$.

If $w{=}(w_j){\in}\N^J$, we define $\|w\|{=}w_1{+}\cdots{+}w_J$ and, for $1{\le}j{\le}J$,  $e_j$ denotes the $j$th unit vector of $\N^J$. It is easily checked that $(X_N(t))$ is an irreducible Markov process on ${\cal S}_N$. Its transition rates are given by 
\begin{itemize}
\item {\sc Transcription of} rRNAs.  For $1{\le}i,j{\le}J$,
  \[
  (f,s,z,(u,r))\longrightarrow
    \begin{cases}
         (f{-}1,s,z,(u,r{+}e_j)) & \alpha_{r,j}\ind{f>0,r_j<C_{r,j}^N},\\
        (0,s,z,(u{-}e_j,r{+}e_j)) & \alpha_{r,j}\ind{u_j=1,r_j<C_{r,j}^N,f{=}0},\\
        (f{+}1,s,z,(u,r{-}e_j)) & \beta_{r,j}\ind{r_j>0,u_k{>}0, \forall 1{\le}k{\le}J},\\
        (0,s,z,(u+e_i,r-e_j)) &\displaystyle \frac{\beta_{r,j}}{ J{-}\|u\|}\ind{r_j>0,u_i=0}.
    \end{cases}
\]
\item {\sc Transcription of} mRNAs.
  \[
  (f,s,z,(u,r))\longrightarrow
    \begin{cases}
         (f{-}1,s,z,(u,r)) & \alpha_{m}f\left(C_m^N{-}\Psi(x)\right),\\
         (f{+}1,s,z,(u,r)) & \beta_{m}\Psi(x).
    \end{cases}
    \]
\item  {\sc Creation/Degradation of} {\sSs}.
  \[
  (f,s,z,(u,r))\longrightarrow
    \begin{cases}
         (f,s,z{+}1,(u,r)) & \beta_6,\\
         (f,s,z{-}1,(u,r)) & \delta_6z.
    \end{cases}
    \]
\item {\sc Sequestration/de-Sequestration of}  Polymerases.
  \[
  (f,s,z,(u,r))\longrightarrow
  \begin{cases}
    (f{-}1,s{+}1,z{-}1,(u,r)) & \lambda f z,\\
    (f+1,s{-}1,z{+}1,(u,r)) & \eta s.
    \end{cases}
  \]
\end{itemize}

\subsection*{A Possible Extension for mRNAs}
We have chosen to consider $C_m^N$ genes of mRNAs with the same parameters $\beta_m$ and $\alpha_m$ for the transcription by polymerases, for simplicity essentially. A generalization could be considered for which the $C_m^N$ types of mRNAs can be split into $K$ sub-groups $({\cal C}^N_{m,k})$ of respective sizes $C^N{\steq{def}}(C^N_{m,k},1{\le}k{\le}K)$ and  with the  parameters $(\beta_{m,k},\alpha_{m,k},1{\le}k{\le}K)$. 

It basically  amounts to say that mRNAs can be partitioned according to the strengths of the affinity of their promoters and of their lengths in terms of nucleotides.  See~\citet{Bremer2008} and~\citet{BioNumbers} for example. A similar assumption for translation of different types of proteins has been done in~\citet{Fromion2}. We denote by $M_{m,j}^N(t)$ the number of polymerases in the elongation phase of an mRNA whose type is in the set ${\cal C}_{m,j}^N$ at time $t$. The state process for this part of the system is 
\[
(M_N(t))\steq{def}(M_{m,k}^N(t)){\in}{\cal S}_m^N, \text{ with } {\cal S}^N_m\steq{def}\left\{x{\in}\prod_{k=1}^{K}\left[0,C_{m,k}^N\right]: |x|\le N\right\}.
\]
Without sequestration and transcription of rRNAs, the model is equivalent to a kind of Ehrenfest urn model, with $K+1$ urns,  for  $1{\le}k{\le}K$, the urn $k$ has a maximal capacity of $C_{m,k}^N$ and the balls inside move to urn $0$ at rate $\beta_{m,k}$. A ball in urn $0$ go to a specific empty place of urn $k$ at rate $\alpha_{m,k}(C_{m,k}^N{-}M_{m,k}^N(t))$ 

\subsection{Orders of Magnitude and Scaling Assumptions}\label{NumbSec}
We now discuss the orders of magnitude of the main parameters of the biological process. 
\begin{itemize}
\item The scaling variable used in our analysis  is  $N$, the total number of polymerases in the cell. It is assumed that this number is constant during the growth phase investigated, this quantity is  quite large, between $2000$ and $10000$ for E. coli, depending of the environment. See~\citet{Bakshi}.
\item The number $J$ of different types of rRNA is  small, of the order of $10$. See~\citet{Bremer2008}.
\item We shall assume that the maximal number of polymerases in transcription of an rRNA of type $j$,  $C_{r,j}^N$, is of the order of $N$, the total number of polymerases. Indeed,  in a steady growth phase a given rRNA gene can accommodate a significant number of polymerases. Recall that the length in nucleotides of an rRNA is large, of the order of $5000$.
\item Similarly, the total number of different types of mRNAs is also of the order of $N$, several thousands, of the order of $3500$ for E. coli.
\end{itemize}
See also~\citet{BioNumbers}, \citet{BioCyc} and~\citet{Neidhardt} for the estimation of the numerical values of these quantities in various contexts.

Due to the modeling assumptions of Section~\ref{ModSec}, we assume that the relations
\begin{equation}\label{Scale}
  \lim_{N\to+\infty} \frac{C_{r,j}^N}{N}=c_{r,j}>0,\  1{\le}j{\le}J, \text{ and }   \lim_{N\to+\infty} \frac{C_m^N}{N}=c_m,
\end{equation}
hold and that, in order to cope with the production of rRNAs during a steady growth phase,  the total number of polymerases $N$ is larger than the total maximal number of polymerases in elongation phase of rRNAs, i.e. that $C_{r,1}^N{+}\cdots{+}C_{r,J}^N$ and, also that there are not too many polymerases for the transcription, i.e.
\[
\max\left(C_m^N,\sum_{j=1}^JC_{r,j}^N\right) < N < C_m^N{+}\sum_{j=1}^JC_{r,j}^N.
\]
In view of~\eqref{Scale}, these assumptions are expressed by the following conditions on the scaled parameters $(c_{r,j})$ and $c_m$,
\begin{equation}\label{SatCond}
\max\left(c_m, \sum_{j=1}^J c_{r,j}\right)< 1 <   \sum_{j=1}^J c_{r,j}{+}c_m.
\end{equation}
We can now introduce the two regimes of interest in our paper. 
\begin{definition}\label{DefPhase}
  \begin{enumerate}
    \item[]
\item The {\em Exponential Phase} is defined by the relation
\begin{equation}\label{ExpCond}
\min_{1{\le}j{\le}J}\frac{\alpha_{r,j}}{\beta_{r,j}}>1.
\end{equation}
The initiation rate $\alpha_{r,j}$ of type $j$  rRNAs is greater than its production rate.
\item  The {\em stationary phase}  is defined by the relation 
\begin{equation}\label{StatCond}
\max_{1{\le}j{\le}J}\frac{\alpha_{r,j}}{\beta_{r,j}}<1. 
\end{equation}
The initiation rate $\alpha_{r,j}$  of type $j$  rRNAs is less than its production rate.
  \end{enumerate}
\end{definition}
It should be noted that Relations~\eqref{ExpCond} and~\eqref{StatCond} are not complementary but this is not a concern for the following reason. If there exists a subset $S$ of $\{1,\ldots,J\}$ such that
\[
\max_{j{\not\in}S}\frac{\alpha_{r,j}}{\beta_{r,j}}<1<\min_{j{\in}S}\frac{\alpha_{r,j}}{\beta_{r,j}},
\]
we will express it as a model for which the rRNAs are defined by the subset $S$ and the remaining nodes $S^c{=}\{1,\ldots,J\}{\setminus}S$ are included in the mRNAs. It can be shown that the addition of a finite number of nodes to the mRNA does not change the first order in $N$ of the number of polymerases in transcription of an mRNA. See Section~\ref{StatSec}.  With this change and if Condition~\eqref{SatCond} holds for this modified system, this is still an exponential phase.
\subsection{An Auxiliary Model}\label{SecAux}
In Sections~\ref{SecSub} and~\ref{SecSuper}, we study a process which can be interpreted as a model similar to $(X_N(t))$  but with only transcription of mRNAs and sequestration by {\sSs} but without rRNAs. The reasons to study this case are two-fold:
\begin{enumerate}
\item Exponential Phase. If Condition~\eqref{ExpCond} holds, as we shall see, ``most'' of the $J{+}C_{r,1}^N{+}C_{r,2}^N{+}{\cdots}{+}C_{r,J}^N$ available places for transcriptions of rRNAs will be occupied by polymerases. Provided that this situation holds on a sufficiently large time scale,  under Condition~\eqref{SatCond}, there are $A_N$ available polymerases for sequestration and transcription of mRNAs, with 
\begin{equation}\label{AN}
A_N\steq{def}N{-}J{-}\sum_{j=1}C_{r,j}^N\sim \gamma N{<} C_m^N.
\end{equation}
The system works as if there were $A_N$ polymerases available for the transcription of mRNAs. With Condition~\eqref{SatCond}, we have $A_N{<}C_m^N$. 
\item Stationary Phase. When Condition~\eqref{StatCond} holds, then, roughly speaking, the total number of polymerases in the elongation phase of rRNAs is $O(1)$, so that this part of the system is in some way negligible. In this case the total number  of polymerases available for transcription of mRNAs is essentially $N$ and thus greater than $C_m^N$ under Condition~\eqref{SatCond}.
\end{enumerate}
The precise definition of  exponential phase, resp. stationary phase, is in Section~\ref{ExpSec}, resp. Section~\ref{StatSec}. 

We denote $(X^0_N(t)){=}(F^0_N(t),S^0_N(t),Z^0_N(t))$ the system defined in Section~\ref{QMatrix} but  without the part of the model for rRNAs. From a state $(f,s,z)$, the transition rates are:  
\begin{equation}\label{Jump1}
    \begin{cases}
      (f{-}1,s,z) & \alpha_{m}\left(C_m^N{-}(N{-}f{-}s)\right)f,\quad \\
      (f{+}1,s,z) & \beta_{m}(N{-}f{-}s),\\
      (f,s,z{+}1) & \beta_6,\\
    \end{cases}
    \begin{cases}
    (f,s,z{-}1) & \delta_6z,\\
    (f{-}1,s{+}1,z{-}1) & \lambda f z,\\
    (f+1,s{-}1,z{+}1) & \eta s.
    \end{cases}
\end{equation}
Using the classical formulation in terms of a martingale problem, see Theorem~(20.6) in Section~IV of~\citet{Rogers1} for example,  the Markov process $(X^0_N(t))$  whose $Q$-matrix is given by Relation~\eqref{Jump1}, as $(F^0_N(t),S^0_N(t),Z^0_N(t))$, the solution  of the SDEs,
\begin{align}
  \diff F^0_N(t)&={\cal P}_1\left(\left(0,\eta S^0_N(t{-}\right), \diff t\right){-}{\cal P}_3\left(\left(0,\lambda F^0_N(t{-})Z^0_N(t{-})\right),\diff t\right)\label{SDE1}\\&{+}{\cal P}_2\left(\left(0,\beta_m\left(N{-}F^0_N(t{-}){-}S^0_N(t{-})\right),\diff t\right)\right)\notag\\& {-}{\cal P}_4\left(\left(0,\alpha_mF^0_N(t{-})\left(C_m^N{-}\left(N{-}F^0_N(t{-}){-}S^0_N(t{-})\right)\right)\right),\diff t\right)\notag \\
\diff S^0_N(t)&={-}{\cal P}_1\left(\left(0,\eta S^0_N(t{-}\right), \diff t\right){+}{\cal P}_3\left(\left(0,\lambda F^0_N(t{-})Z^0_N(t{-})\right),\diff t\right)\label{SDE2}\\
\diff   Z^0_N(t)&={\cal P}_5\left(\left(0,\beta_6\right),\diff t\right)
{-}{\cal P}_6\left(\left(0,\delta_6Z^0_N(t-)\right),\diff t\right)\label{SDE3}\\
&{+}{\cal P}_1\left(\left(0,\eta S^0_N(t{-}\right), \diff t\right) {-}{\cal P}_3\left(\left(0,\lambda F^0_N(t{-}) Z^0_N(t{-})\right),\diff t\right),\notag
\end{align}
with the convenient initial conditions, where ${\cal P}_i$, $i{\in}\{1,2,3,4\}$ are independent Poisson processes on $\R_+^2$ with intensity $\diff s{\otimes}\diff t$. 

We will study two regimes of this stochastic model:
\begin{itemize}
\item Sub-critical case, when $c_m{>}1$, i.e. $N{<}C_m^N$ for $N$ sufficiently large.
\item Super-critical case, when $c_m{<}1$.
\end{itemize}
As it will be seen these two regimes are respectively associated to the exponential and stationary phases.

\subsection*{Notations}
We define a common filtration common to all our processes, as follows, for $t{\ge}0$
\begin{equation}\label{Filt}
{\cal F}_t=\sigma\left({\cal P}_i(A{\times}[0,s]): A{\in}{\cal B}(\R_+), i{\in}\{1,\cdots,6\}, s\le t\right).
\end{equation}
From now on, all notions of stopping time, adapted process,  martingale, refer to this (completed) filtration. A \cadlag process is an adapted process such that with probability one, it is right-continuous process with left limits at any positive real number. 

If $H$ is a locally compact metric space, we denote by ${\cal C}_c(H)$ the set of continuous functions with compact support on $H$. It is endowed with the topology of the uniform norm. The set ${\cal P}(H)$ is the space of Borelian probability distributions on~$H$.
\section{Sub-critical Case}\label{SecSub}
It is assumed throughout this section that $c_m{>}1$ holds, the total number of possible sites for transcription of mRNAs is much larger than the total number of polymerases.

\begin{definition}[Occupation measure of $(F^0_N(t))$]\label{OccMu}
For $g{\in}{\cal C}_{c}\left(\R_+{\times}\N\right)$
\begin{equation}\label{EmpSub}
    \croc{\mu_N,g}\steq{def} \int_0^{+\infty} g\left(u,{F}^0_N(u)\right)\diff u.
\end{equation}
\end{definition}

We start with a technical result on a birth and death process. 

\begin{lemma}\label{lemMM}
For $\kappa_i{>}0$ and $\kappa_o{>}0$,  let $(Y(t))$ be a birth and death process on $\N$ whose $Q$-matrix is given by
\[
  q(x,x{+}1)=\kappa_i \text{ and }   q(x,x{-}1)=\kappa_o x,\quad x{\in}\N,
\]
\begin{enumerate}
\item if $Y(0){=}N$ and 
\[
H_Y^N\steq{def}\inf\{t{>}0: Y(t){=}0\},
\]
then $(H_Y^N/\ln N)$ is converging in distribution to a constant.
\item if $Y(0){=}0$, then for any  $\delta{>}0$, the convergence in distribution 
\[
\lim_{N\to+\infty}\left(\frac{Y\left(N^{\delta}t\right)}{\ln(N)^2}\right)=0
\]
holds.
\end{enumerate}
\end{lemma}
The process  $(Y(t))$ can be thought as a kind of discrete Ornstein-Uhlenbeck process on $\N$. In a queueing context, this is the process of the number of jobs of an $M/M/\infty$ queue. See Chapter~6 of~\citet{Robert} for example. Its invariant distribution is Poisson with parameter $\kappa_i/\kappa_o$. 
\begin{proof}
The first assertion comes directly from Proposition~6.8 of~\cite{Robert}. If $Y(0){=}0$ and, for $p{\ge}1$,
  \[
  T_p=\inf\{t{>}0:Y(t)>p\},
  \]
  Proposition~6.10~\cite{Robert} gives that, if $\rho{\steq{def}}\kappa_i/\kappa_o$, the sequence
  \[
  \left(\rho^p\frac{T_p}{p!}\right)
  \]
  is converging in distribution to an exponential distribution. In particular for any $K{>}0$,
  \[
  \lim_{N\to+\infty}\P\left(T_{\ln(N)^2}<KN^\delta\right)=0,
  \]
  since,   by Stirling's Formula, 
  \[
  \lim_{N\to+\infty} \frac{(\ln(N)^2)!}{\rho^{\ln(N)^2}N^\delta}={+}\infty.
  \]
The lemma is proved.
\end{proof}

We begin with a lemma showing that the initial state of $({X}^0_N(t))$ can be taken with few free polymerases. 
\begin{lemma}\label{LemMIHit}
  If $c_m{>}1$ and  $({F}^0_N(0),{S}^0_N(0),{Z}^0_N(0)){=}(f_N,s_N,z_N)$, such that
  \[
  \lim \frac{1}{N}(f_N,s_N,z_N)=((f_0,s_0,z_0){\in}\R_+^3,
  \]
with $f_0{+}s_0{<}1$,  and, if
    \[
  {\tau}^0_N\steq{def}\inf\left\{t:{F}^0_N(t)=0\right\},
  \]
  then the sequence $({\tau}^0_N/(\ln N)^2)$ is converging in distribution to $0$. 
\end{lemma}
The Condition $f_0{+}s_0{<}1$ is to take into account  the fact that ${F}^0_N(0){+}{S}^0_N(0){\le}N$. 
\begin{proof}
Because of the assumption  $c_m{>}1$,  for $N$ sufficiently large, there exists $\eps_0{>}0$ such that $C_m^N{-}N{>}\eps_0N$, the relations~\eqref{Jump1} for the transition rates show that one can construct a coupling $({F}^0_N(t),Y(t))$ such that $Y(0){=}{F}^0_N(0)$ and the relation
\begin{equation}\label{UpMMI}
  {F}^0_N(t)\le Y(Nt),\qquad \forall t{\ge}0,
\end{equation}
  holds almost surely for all $t{\ge}0$,   where $(Y(t))$ is a process as defined in Lemma~\ref{lemMM} with $\kappa_i{=}\eta{+}\beta_m$ and $\kappa_o{=}\alpha_m \eps_0$. We conclude the proof by using Lemma~\ref{lemMM}.
\end{proof}
We can now state the main result of this section.  It shows that for the asymptotic system, when $N$ is large,  all polymerases are eventually in the transcription phase of mRNAs, i.e. the fraction of sequestered polymerases is close to~$0$. A sketch of the proof is given in Section~\ref{proofpropsub} of the Appendix.  
\begin{proposition}[Starting from a Congested State]\label{propsub}
Under the condition $c_m{>}1$ and  if the initial state is $({F}^0_N(0),{S}^0_N(0),{Z}^0_N(0)){=}(0,s_N,z_N)$ and 
 \[
  \lim \frac{1}{N}(s_N,z_N)=(s_0,z_0)\in\R_+^2,
  \]
  such that $s_0{+}z_0{<}1$,   then, for the convergence in distribution
  \[
  \lim_{N\to+\infty} \left(\frac{{S}^0_N(t)}{N},\frac{{Z}^0_N(t)}{N}\right)=(s(t),z(t)),
  \]
  where $(s(t),z(t))$ is the unique solution of the system of ODEs,
\[
  \dot{s}(t)={-}\eta  s(t){+}\lambda z(t) \frac{\beta_m{-}\left(\beta_m{-}\eta\right)s(t)}{\alpha_m(c_m{-}1{+}s(t))+\lambda z(t)},
  \quad 
\dot{s}(t){+}\dot{z}(t)={-}\delta_6 z(t),
\]
with $(s(0),z(0))=(s_0,z_0)$.
\end{proposition}
It is not difficult to see that the function $(s(t),z(t))$ is converging to $(0,0)$ at infinity.

To study the asymptotic behavior of the model in the exponential regime, we investigate the occupation measure associated to free polymerases when the initial state is ``small''. In this case, contrary to the last proposition, the processes $(S^0_N(t),Z^0_N(t))$ should be ``slow'', i.e. their transition rate are of the order of $O(1)$, only $(F^0_N(t))$ is ``fast''. 

\begin{proposition}[Fixed Initial Point]\label{propsub2}
  Under the condition $c_m{>}1$ and the initial state is such that  $({F}^0_N(0),{S}^0_N(0),{Z}^0_N(0)){=}(f_0,s_0,z_0){\in}\N^3$,
  then, for the convergence in distribution
    \[
    \lim_{N\to+\infty} \croc{\mu_N,g}= \int_0^{+\infty} \E\left(g\left(u,{\cal N}_1\left(0,\rho_m\right)\right)\right)\diff u,
  \]
  for any $g{\in}{\cal C}_{c}\left(\R_+{\times}\N\right)$, where $\rho_m{=}{\beta_m}/{(\alpha_m(c_m-1))}$, $\mu_N$ is the occupation measure defined by Relation~\eqref{EmpSub}, and ${\cal N}_1$ is a Poisson process with rate $1$.

The sequence of processes $(S^0_N(t),Z^0_N(t))$ converges in distribution for the Skorohod topology to a jump process  $(Y(t))$ on $\N^2$ whose transition rates are given by
  \[
  (s,z)\longrightarrow (s,z){+}
  \begin{cases}
    (1,{-}1)& \lambda\rho_m z,\\
    ({-}1,1)& \eta s,
  \end{cases}
  \quad 
  \begin{cases}
    (0,1)& \beta_6,\\
    (0,{-}1)& \delta_6z.
  \end{cases}
  \]
\end{proposition}
See Section~\ref{proofpropsub2} of the appendix. 
\section{Super-critical Case}\label{SecSuper}
In this section we study the auxiliary process under the condition $c_m{<}1$, so that  $C_m^N{<}N$ for $N$ sufficiently large. In this case the places for transcription of mRNAs are likely to be saturated quickly. Consequently, there should remain many free polymerases and the sequestration mechanism  has to play a role.

If there are no {\sSs} initially, since the creation of {\sSs} is constant,  the  sequestration of a significant fraction of these polymerases will occur after a duration of time at least of the order of $N$.  In this case, when a {\sS} is created, it  is  right away paired with a free polymerase and will paired again and again that after the successive steps of sequestration/desequestration,  as long as the number of free polymerases is sufficiently ``large''. The sequestration occurs always before a possible degradation of the {\sSs} takes place. The precise result is in fact more subtle than that. It will be shown that, on the fast time scale $t{\mapsto}N t$, the sequestration of polymerases increases but, due to the degradation of {\sSs} there will remain a positive fraction of free polymerases. 

The goal of this section is of proving an averaging principle for the process $(F_N^0(t),Z_N^0(t))$. A coupling and a technical lemma are presented in Section~\ref{Subtechlem}, tightness properties of occupations measures are proved in Section~\ref{Subtight}, finally the main convergence results are proved  in Section{Subavpr}. 

\begin{definition}
  For $N{>}0$ and $t{\ge}0$, we define
  \[
  G^0_N(t)\steq{def}C_m^N{-}\left(N{-}F^0_N(t){-}S^0_N(t)\right),
  \]
  the number of ``empty'' places for transcription of mRNAs at time $t$. 

 The scaled process is defined by  
\begin{equation}\label{Fbar}
  \left(\overline{X}^0_N(t)\right)=\left(\overline{F}^0_N(t),G^0_N(Nt),Z^0_N(Nt)\right) \text{ with }
  \left(\overline{F}^0_N(t)\right){\steq{def}}\left(\frac{F^0_N(N t)}{N}\right). 
\end{equation}
If $g$ is non-negative Borelian function on $\R_+^2{\times}\N^2$, we define the occupation measure
\begin{equation}\label{OccMeasDef}
  \croc{\overline{\Lambda}^0_N,g}\steq{def} \int_{\R_+}g\left(s,\overline{X}^0_N(s)\right)\diff s.
\end{equation}
\end{definition}
Not that the, a priori, slow process $(\overline{F}_N(t))$ is also included in the definition of the occupation measure $\overline{\Lambda}^0_N$. The reason is that the proof of the tightness of $(\overline{F}_N(t))$ (for the topology of the uniform norm on \cadlag functions) is not clear.  Due to the fast time scale, the proof that the martingale component of $(\overline{F}_N(t))$ vanishes does not seem to be straightforward. 

The following initial conditions will be assumed,
\begin{equation}\label{InitOp}
\lim_{N\to+\infty} \frac{F^0_N(0)}{N}{=}\overline{f}_0{\in}(0,1{-}c_m),\, G^0_N(0){=}m_0,\, \text{ and }  Z^0_N(0){=}z_0,
\end{equation}
with $m_0$ $z_0{\in}\N$.  A fraction $\overline{f}_0$ of the polymerases are initially free and there are $z_0$ {\sSs} and  $C_m^N{-}m$ polymerases in the transcription phase of mRNAs and the number of sequestered polymerases $S^0_N(0)$ is therefore such that
\[
\lim_{N\to+\infty} \frac{S^0_N(0)}{N}{=}1{-}\overline{f}_0{-}c_m.
\]
As it will be seen in Section~\ref{Subtechlem} there is no loss of generality to consider these initial conditions. 

Before proving the convergence of the sequence of processes $(\overline{X}_N^0(t))$, we analyze the convergence of a ``stopped'' version of it. In several technical arguments we will need that the fraction of free polymerases is not too small. A second step is of showing that, essentially,  the stopped process does not differ from the original process. 
\begin{definition}
For $a{>}0$, the stopping time $\tau_N(a)$ is defined by 
\begin{equation}\label{tau}
\tau_N(a)\steq{def}\inf\left\{t{>}0: F^0_N(N t) \le a N\right\},
\end{equation}
and
\begin{enumerate}
\item if $(W(t))$ is a \cadlag process, we denote $(W_N^a(t)){=}(W(N(t{\wedge}\tau_N(a))))$;
\item The ``stopped'' occupation measure $\overline{\Lambda}^{0,a}_N$ is  defined by, if $g$ is non-negative Borelian function on $\R_+{\times}\N{\times}\R_+$,
\[
  \croc{\overline{\Lambda}_N^{0,a},g}\steq{def} \int_{0}^{\tau_N(a)}g\left(s,\overline{X}^0_N(s)\right)\diff s.
\]
\end{enumerate}
\end{definition}
With a slight abuse, the notation $(\overline{F}_N^a(t)){=}(F^0_N(N(t{\wedge}\tau_N(a)))/N)$ will be used in the following. 

\subsection{Technical Lemmas}\label{Subtechlem}
The two processes  $(G^0_N(t)$ and $(Z^0_N(t))$ are in fact in a neighborhood of $0$ quickly. They will be the fast processes (on the timescale $t{\mapsto}Nt$) of our averaging principle. In state $F^0_N{=}f$, $G^0_N{=}g$ and $Z^0_N{=}z$, $f$, $g$, $z{\in}\N$, the jump rates of the process $(G^0_N(t)$ and $(Z^0_N(t))$ are respectively
\[
\begin{cases}
{+}1,&  \beta_m(C_m^N{-}g),\\
{-}1,& \alpha_m f g,
  \end{cases}
\quad  \text{ and } \quad 
 \begin{cases}
{+}1,&  \left(\beta_6{+}\eta (N{-}f-C_m^N+g)\right),\\
{-}1 ,& \left(\lambda  f{+}\delta_6\right)z.
  \end{cases}
  \]
  If $\eta_0{>}\eta$ and $\eta_1{>}\beta_mc_m$, and $N$ sufficiently large, up to time $\tau_N(a)$, a simple coupling shows that there exist  independent processes $(Y_G(t))$ and $(Y_Z(t))$ such that
\begin{equation}\label{CoupGZ}
     G^0_N(Nt)\le Y_G(N^2t) \text{ and } Z^0_N(Nt)\le Y_Z(N^2t),
\end{equation}
holds for all $t{\in}(0,\tau_N(a))$. The process  $(Y_G(t))$, resp. the process $(Y_Z(t))$, is as in Lemma~\ref{lemMM} with $\kappa_{i,G}{=}\eta_1$ and  $\kappa_{o,G}{=}\alpha_m a$ (resp. $\kappa_{i,Z}{=}\eta_0$ and  $\kappa_{o,Z}{=}\lambda a$), and  $Y_G(0){=}G^0_N(0)$,  resp. $Y_Z(0){=}Z^0_N(0)$. It is not difficult, using again Lemma~\ref{lemMM}, as in Section~\ref{SecSub}, that the hitting time of $(0,0)$ by $(Y_G(t),Y_Z(t))$ is of the order of $\ln N$ so that Condition~\eqref{InitOp} for the initial state can be assumed. 

\begin{lemma}\label{lem1Op}
Under Conditions~\eqref{InitOp} and  $c_m{<}1$,  and if   $a{\in}(0,\overline{f}_0)$, then
  \[
  \lim_{N\to+\infty} \P\left(\tau_N(a){<} t_0^a\right)=0,
  \]
with $t_0^a{=}(\overline{f}_0{-}a)/\beta_6$, and the relation
\[
\lim_{N\to+\infty} \left(\frac{G^{0,a}_N(Nt)}{\ln(N)^2},\frac{Z^{0,a}_N(Nt)}{\ln(N)^2}\right)=(0).
\]
 holds for the convergence in distribution. 
\end{lemma}
In the following, we will use the notation $t_0^a$, where $a{\in}(0,\overline{f}_0)$ is fixed. 
\begin{proof}
The first relation is clear since, for $x{\le}1$ and $t{>}0$, on the event $\{\overline{F}^0_N(t){<}x\}$ there are at least $F^0_N(0){-}\lfloor Nx\rfloor {-}z_0$ new {\sSs} which have been created up to time~$t$.
The rest of the proof follows from the coupling with $(Y_G(t),Y_Z(t))$, Relation~\eqref{CoupGZ},  and Lemma~\ref{lemMM}.
\end{proof}

\subsection{Tightness Properties}\label{Subtight}
\begin{proposition}\label{PropLa}
Under Conditions~\eqref{InitOp} and $c_m{<}1$, the sequence of measure-valued processes $(\overline{\Lambda}_N^{0,a})$ on the state space $[0,t_0^a){\times}\R_+{\times}\N^2$ is tight for the convergence in distribution and any limiting point $\overline{\Lambda}_\infty^{0,a}$ can be expressed as,
\begin{equation}\label{SAPpi}
\croc{\overline{\Lambda}_\infty^{0,a},f}=\int_{[0,t_0^a){\times}\R_+{\times}\N^2} f\left(s,x,p\right)\pi_s^a(\diff x,\diff p)\diff s,
\end{equation}
for any function $f$ with compact support on $[0,t_0^a){\times}\R_+{\times}\N^2$,
where  $t_0^a{=}(\overline{f}_0{-}a)/\beta_6$ and  $(\pi_s^a)$ is an optional process with values in ${\cal P}(\R_+{\times}\N^2)$.
\end{proposition}
 For an introduction to the convergence in distribution of measure-valued processes, see~\citet{Dawson}. The optional property is just used to have convenient measurability properties to define random variables as integrals with respect to $(\pi_s^a,s{>}0)$. See Section VI.4 of~\citet{Rogers2}.

\begin{proof}
Note that, for $K{>}0$ and $t{<}t_0^a$, since
\[
 \int_0^t\ind{Z^0_N(Ns)\ge K}\diff s= \int_0^t\ind{Z^{0,a}_N(s)\ge K}\diff s
 \]
 holds on the event $\{\tau_N(a){\ge}t_0^a\}$, then
\begin{multline*}
\E\left( \overline{\Lambda}^{0,a}_N([0,t_0^a]{\times}[0,1]{\times}\N{\times}[K,{+}\infty])\right)\\\le\E\left(\ind{\tau_N(a)>t_0^a}\int_0^{t_0^a}\ind{Z^{0,a}_N(s)\ge K}\diff s\right) {+}t_0^a\P\left(\tau_N(a){<}t_0^a\right),
\end{multline*}
and, with Relation~\eqref{CoupGZ} and Lemma~\ref{lem1Op}, we have
\begin{multline*}
\E\left(\ind{\tau_N(a)>t_0^a}\int_0^{t_0^a}\ind{Z^{0,a}_N(s)\ge K}\diff s\right)\\  \le \int_0^{t_0^a}\P(Y_Z(N^2s)\ge K)\diff s=\frac{1}{N^2}\int_0^{N^2t_0^a}\P(Y_Z(s)\ge K)\diff s,
\end{multline*}
since the Markov process $(Y_Z(t))$ converges in distribution to a Poisson distribution with parameter~$\kappa_{i,Z}/\kappa_{o,Z}$, the ergodic theorem for Markov processes and   Lemma~\ref{lem1Op}  give therefore the inequality
\[
\limsup_{N\to+\infty}   \E\left( \overline{\Lambda}^{0,a}_N([0,t_0^a]{\times}[0,1]{\times}\N{\times}[K,{+}\infty])\right)\le
t_0^a \P({\cal N}_1(0,\eta/(\lambda a)){\ge}K),
\]
where ${\cal N}_1$ is a   Poisson process on $\R_+$ with rate $1$. One can choose $K$ sufficiently large such that
$\E\left( \overline{\Lambda}^{0,a}_N([0,t_0^a]{\times}[0,1]{\times}\N{\times}[K,{+}\infty])\right)$ is arbitrarily small uniformly in $N$. Similarly, by replacing $(Z^0_N,Y_Z)$ by $(G^0_N,Y_G)$  the same property can be proved for $\E\left( \overline{\Lambda}^{0,a}_N([0,t_0^a]{\times}[0,1]{\times}[K,{+}\infty]{\times}\N\right)$ for $K$ and $N$ sufficiently large.
For any $\eps{>}0$, there exists some $K_0$ such that
\[
\sup_{N} \E\left( \overline{\Lambda}^{0,a}_N([0,t_0^a]{\times}[0,1]{\times}[0,K_0]^2\right)\ge (1{-}\eps)t_0^a.
\]
Lemma~1.3 of~\citet{Kurtz} shows that the sequence $(\overline{\Lambda}^{0,a}_N)$ is  tight, and Lemma~1.4 of the same reference gives the representation~\eqref{SAPpi}.
\end{proof}

Proposition~\ref{PropLa} has established tightness properties   $(\overline{\Lambda}_N^{0,a})$. The following simple lemma extends this result in terms of the convergence of stochastic processes. It will be used  repeatedly, in particular to identify the possible limits of  $(\overline{\Lambda}_N^{0,a})$. See~\citet{Dawson} for example.
\begin{lemma}\label{ContLem}
Under Conditions~\eqref{InitOp}  and $c_m{<}1$, if $(\overline{\Lambda}_{N_k}^{0,a})$ is a subsequence converging to $\overline{\Lambda}^{0,a}_\infty$ satisfying Relation~\eqref{SAPpi}, then for any $g{\in}{\cal C}_c(\R_+{\times}\N^2)$, for the convergence in distribution of processes associated to the uniform norm, 
\[
\lim_{k\to+\infty} \left( \int_0^t g\left(\overline{X}^0_{N_k}(s)\right)\diff s\right)=
 \left( \int_0^t\int_{\R_+{\times}\N^2} g\left(x,p\right)\pi_s^a(\diff x,\diff p)\diff s \right). 
 \]
\end{lemma}
\begin{proof}
  The tightness of the sequence of stochastic processes is obtained by the use of the criterion of the modulus of continuity. See Theorem~7.3 of~\citet{Billingsley}. The identification of the limit is a straightforward  consequence of the convergence of  $(\overline{\Lambda}_{N_k}^{0,a})$
\end{proof}
If we divide by $N^2$ Relation~\eqref{SDEf} of the appendix, we get that, on the event $\{\tau_N(a){>}t\}$, the relation
\begin{align}
\frac{1}{N^2}&f\left(\overline{X}^0_N(t)\right)=\frac{1}{N^2}f\left(\overline{X}^0_N(0)\right){+}\frac{M_{f,N}(t)}{N^2}\label{SDEf2}\\
&+\lambda \int_0^{t} \nabla_{-\frac{e_1}{N}-e_3}(f)\left(\overline{X}^0_N(s)\right) \overline{F}^0_N(Ns)Z^0_N(Ns)\diff s\notag\\
&  {+}\eta \int_0^{t}\nabla_{\frac{e_1}{N}{+}e_3}(f)\left(\overline{X}^0_N(s)\right)\left(1{-}\frac{C_m^N}{N}{+}\frac{G^0_N(Ns)}{N}{-}\overline{F}^0_N(Ns)\right)\diff s\notag\\
&+\alpha_m \int_0^{t} \nabla_{-\frac{e_1}{N}-e_2}(f)\left(\overline{X}^0_N(s)\right)G^0_N(Ns)\overline{F}^0_N(Ns)\diff s\notag\\
&+\beta_m  \int_0^{t} \nabla_{\frac{e_1}{N}{+}e_2}(f)\left(\overline{X}^0_N(s)\right)\left(\frac{C_m^N}{N}{-}\frac{G^0_N(Ns)}{N}\right)\diff s\notag\\
&+\frac{\beta_6}{N} \int_0^{t} \nabla_{e_3}(f)\left(\overline{X}^0_N(s)\right)\diff s+\frac{\delta_6}{N} \int_0^{t} \nabla_{-e_3}(f)\left(\overline{X}^0_N(s)\right)Z^0_N(Ns)\diff s\notag
\end{align}
holds. Recall that, for $i{\in}\{1,2,3\}$, $e_i$ is the $i$th unit vector of $\R^3$. 
\begin{lemma}\label{LemMart}
Under Conditions~\eqref{InitOp} and $c_m{<}1$,  if $f$ is a continuous bounded function on $\R_+{\times}\N$, then the martingale $(M_{f,N}(t)/N^2, t{<}t_0^a)$ of Relation~\eqref{SDEf2} converges in distribution to $0$. 
\end{lemma}
\begin{proof}
We take care of one of the six terms of $(\croc{M_{f,N}/N^2}(t))$ of Relation~\eqref{CrocMf}, the arguments are similar for the others, even easier. 
\[
A_{1,N}(t)\steq{def}\frac{\lambda}{N^2} \int_0^{t} \left[\nabla_{-\frac{e_1}{N}-e_3}(f)\left(\overline{X}^0_N(s)\right)\right]^2\overline{F}^0_N(s)Z^0_N(Ns)\diff s
\]
We note that for $t{\ge} 0$, $0{\le}{Z}^0_N(t){\le}N{+}{\cal P}_5((0,\beta_6){\times}(0,t])$. Consequently, Doob's Inequality shows the convergence of $(M_{f,N}(t)/N^2)$ to $0$. The lemma is proved.

\end{proof}
\begin{proposition}\label{prop12}
Under Conditions~\eqref{InitOp}  and $c_m{<}1$, and if $\overline{\Lambda}_\infty^{0,a}$ is a limiting point of $\overline{\Lambda}_n^{0,a}$ with the representation~\eqref{SAPpi} of Proposition~\ref{PropLa}, then, if $\pi^{1,a}_t{=}\pi_t^{0,a}(\cdot,\N^2)$, for any $t{<}t_0^a$ and any continuous function $g$ on  $\R_+{\times}\N^2$ we have
\begin{multline}\label{SAPpi2}
\int_0^t \int g(x,p)\pi_s^a(\diff x,\diff p)\diff s\\ =\int_0^t \int_{\R_+} \E\left[g\left(x,{\cal N}_1\left(\left[0,\rho_m\frac{c_m}{x}\right]\right), {\cal N}_2\left(\left[0,\rho_1\frac{1{-}c_m{-}x}{x}\right]\right)\right)\right]\pi_s^{1,a}(\diff x)\diff s,
\end{multline}
where ${\cal N}_1$ and  ${\cal N}_2$  are two independent   Poisson processes on $\R_+$ with rate $1$ and
\begin{equation}\label{rho}
\rho_1{=}\frac{\eta}{\lambda} \text{ and }\rho_m{=}\frac{\beta_m}{\alpha_m}.
\end{equation}
\end{proposition}
Relation~\eqref{SAPpi2} states that, for almost all  $t{<}t_0^a$,  $\pi_t$ conditioned on the first coordinate $x$ is a product of two Poisson distributions with respective parameters $\rho_mc_m/x$ and $\rho_1(1{-}c_m{-}x)/x$. 
\begin{proof}
Let $(\overline{\Lambda}^0_{N_k})$ be a subsequence of $(\overline{\Lambda}^0_N)$ converging to some $\overline{\Lambda}^0_\infty$ of the form~\eqref{SAPpi}. 
By letting $k$ go to infinity in Relation~\eqref{SDEf2}, with Lemmas~\ref{lem1Op}, \ref{ContLem} and~\ref{LemMart},  we obtain that there exists an event ${\cal E}_1$ with $\P({\cal E}_1){=}1$  on which the relation
\begin{multline*}
\int_0^{t}\int_{\R_+{\times}\N^2} \left(\eta\left(1{-}c_m{-}x\right)\nabla_{e_3}(f)(x,p){+} \lambda xp_2\nabla_{-e_3}(f)(x,p) \right)\pi_s^{a}(\diff x,\diff p)\diff s\\
+ \int_0^{t}\int_{\R_+{\times}\N^2} \left(\beta_mc_m  \nabla_{e_2}(f)(x,p){+}\alpha_m\nabla_{-e_2}(f)(x,p) p_1x\right)\pi_s^{a}(\diff x,\diff p)\diff s=0,
\end{multline*}
holds for all $t{\le}T$ and for all  functions $f{\in}{\cal C}_c(\R_+{\times}\N^2)$, by using the separability property of this space of functions for the uniform norm.   If $f(x,p){=}f_1(x)f_2(p)$, this relation can be rewritten as
\[
\int_0^t\int_{\R_+{\times}\N^2}  f_1(x)\Omega[x](f_2)(p)\pi_s^{a}(\diff x,\diff p)\diff s=0
\]
where, for $h:{\N^2}{\to}\R_+$ and $p{=}(p_1,p_2){\in}\N^2$, 
\begin{multline*}
\Omega[x](h)(p)=\beta_mc_m \nabla_{e_1}(h)(p){+}\alpha_m p_1x\nabla_{-e_1}(h)(p) \\
{+}\eta(1{-}c_m{-}x)\nabla_{e_2}(h)(p) {+}\lambda xp_2  \nabla_{-e_2}(h)(p).
\end{multline*}
$\Omega[x]$ is the jump matrix of two independent birth and death processes $(Y_1(t))$ and $(Y_2(t))$ as  in Lemma~\ref{lemMM}
with  parameters $\kappa_i{=}\beta_mc_m$, $\kappa_o{=}\alpha_m x$ for $(Y_1(t))$   and $\kappa_i{=}\eta(1{-}c_m{-}x)$, $\kappa_o{=}\lambda x$  for $(Y_2(t))$. 

Consequently, for almost all $t{\le}T$, the relation
\[
\int_{\R_+{\times}\N^2}  f_1(x)\Omega[x](f_2)(p)\pi_t^a\diff x,\diff p)=0
\]
holds. Hence, if $\widetilde{\pi}^a_t(\cdot|x)$ is the conditional probability on $\N^2$ of $\pi^a_t(\diff x,\diff p)$ given $x$, we have
\[
\int_{\R_+}  f_1(x)\int_{\N^2}\Omega[x](f_2)(p)\widetilde{\pi}^a_t(\diff p|x)\pi^{1,a}(\diff x)=0,
\]
we deduce that the relation
\[
\int_{\N^2}\Omega[x](f_2)(p)\widetilde{\pi}^a_t(\diff p|x)=0
\]
holds $\pi_t^{1,a}(\diff x)$ almost surely,  for all functions $f_2$ with finite support on $\N^2$. Consequently,  $\pi^{1,a}(\diff x)$ almost surely, $\widetilde{\pi}^a_t(\diff p|x)$ is the invariant distribution associated to the $Q$-matrix $\Omega[x]$.   The proposition is proved. 
\end{proof}
We fix $(N_k)$  an increasing sequence such the sequence $(\overline{\Lambda}_{N_k}^{0,a})$ is converging in distribution to the law of  $\overline{\Lambda}_\infty^{0,a}$ with a representation given by Relations~\eqref{SAPpi} and~\eqref{SAPpi2}. 

\subsection{Averaging Principle}\label{Subavpr}
We define, for $t{\ge}0$, 
\[
\widetilde{Z}^0_N(t)=S^0_N(t){+}Z^0_N(t),
\]
$\widetilde{Z}^0_N(t)$ is in fact the total number of {\sSs} (free or paired)  of the system at time~$t$.  Using the SDEs~\eqref{SDE2} and~\eqref{SDE3}, we have
\begin{equation}\label{eqZb}
\frac{\widetilde{Z}^0_N(Nt)}{N}=M_{Z,N}(t){+}\frac{\widetilde{Z}^0_N(0)}{N}{+}\beta_6 t{-}\delta_6\int_0^tZ^0_N(Ns)\diff s,
\end{equation}
where $(M_{Z,N}(t))$ is a local martingale whose previsible increasing process is given by
\begin{equation}\label{eqZMcroc}
\left(\croc{M_{Z,N}}(t)\right)=\left(\frac{1}{N}\left(\beta_6 t{+}\delta_6\int_0^tZ^0_N(Ns)\diff s\right)\right).
\end{equation}

\begin{proposition}\label{CVOcc1}
Under Conditions~\eqref{InitOp}  and $c_m{<}1$, for the convergence in distribution
\[
\lim_{k\to+\infty}  \left(\int_0^t Z_{N_k}(N_ks)\diff s,t<t_0^a\right)=\left(\rho_1\int_0^t\int_{\mathbb{R}_+}\frac{1{-}c_m{-}x}{x}\pi_s^1(\diff x)\diff s, t< t_0^a\right),
\]
with $t_0^a{=}(\overline{f}_0{-}a)/\beta_6$ and $\rho_1{=}\eta/\lambda$.  Furthermore, $(M_{Z,N}(t),t{<}t_0^a)$ is converging to $0$.
\end{proposition}
\begin{proof}
 The convergence of the sequence of stochastic processes  $(M_{Z,N}(t),t{<}t_0^a)$ to $0$ is a consequence of Relations~\eqref{eqZb} and~\eqref{eqZMcroc}, and of Doob's Inequality. 
For $0{\le}s{\le}t$, the coupling~\eqref{CoupGZ} and Cauchy-Schwartz' Inequality give
\begin{multline*}
\E\left(\left(\ind{\tau_N(a){>}t}\int_s^t Z_{N_k}(Ns)\diff s \right)^2\right)\leq (t-s)\E\left(\ind{\tau_N(a){>}t}\int_s^t Z_{N_k}^a(s)^2\diff s \right)\\ \le (t{-}s)\E\left(\int_s^t Y_Z(Ns)^2\diff s \right)\le (t{-}s)^{2}\sup_{u{\ge}0}\E\left(Y_Z(u)^2\right).
\end{multline*}
We now use the Kolmogorov-\v{C}entsov's criterion,  see Theorem~2.8 and Problem~4.11, page~64 of~\citet{Karatzas} and Lemma~\ref{lem1Op} to show  that the sequence of stochastic processes
\[
\left(\int_0^t Z_{N_k}(Ns)\diff s, t<t_0^a\right)
\]
is tight for the convergence in distribution.

Lemma~\ref{ContLem} gives the convergence in distribution
\begin{multline*}
\lim_{k\to+\infty} \left(\int_0^t Z_{N_k}(N_ks){\wedge}K \diff s,t<t_0^a\right) \\=\left(\int_0^t\int_{\mathbb{R}_+}\E\left({\cal N}_1\left(0,\rho_1\frac{1{-}c_m{-}x}{x}\right){\wedge}K\right)\,\pi_s^1(\diff x)\diff s,t<t_0^a\right). 
\end{multline*}
By using again Relation~\eqref{CoupGZ},   we have
\[
\E\left(\int_0^{t_0^a} Z_{N_k}(N_ks)\ind{Z_{N_k}(N_ks)\ge K} \diff s\right)
\leq \E\left(\int_0^{t_0^a} Y_Z({N_k}^2s)\ind{Y_Z(N_k^2s)\ge K} \diff s\right)
\]
and the convergence in distribution of $(Y_Z(t))$, as $t$ goes to infinity, to  $Y_Z(\infty)$ a random variable  with a Poisson distribution with parameter $\rho_Z{=}\kappa_{i,Z}/\kappa_{o,Z}$ give 
\[
\lim_{k\to+\infty} \E\left(\int_0^{t_0^a} Y_Z(N_k^2s)\ind{Y_Z(N_k^2s)\ge K} \diff s\right)
=t_0^a\E\left(Y_Z(\infty)\ind{Y_Z(\infty){\ge}K}\right).
\]
 It is then easy to obtain the first convergence by letting $K$ go to infinity.

The proposition is proved. 
\end{proof}

Relation~\eqref{eqZb} therefore shows that, on the time interval $I_a{=}[0,t_0^a)$, the sequence of processes
\[
\left(\frac{\widetilde{Z}^0_{N_k}(N_kt)}{N_k}\right)
\]
is converging in distribution. Since
\[
\left(\frac{\widetilde{Z}^0_{N_k}(N_kt)}{N_k}\right)=\left(1{-}\frac{F_{N_k}^0(N_kt)}{N_k}{-}\frac{C_m^{N_k}}{N_k}{+}\frac{G_{N_k}^0(N_kt)}{N_k}{+}\frac{Z_{N_k}^0(N_kt)}{N_k}\right)
\]
with  Lemma~\ref{lem1Op}, we therefore obtain that the  sequence of processes $(F_{N_k}(N_kt)/N_k)$ is converging in distribution to some process $(\overline{f}(t))$ on $I_a$. In particular, for $t{<}t_0^a$ and $g{\in}{\cal C}_c(\R_+)$,  we have
\[
\int_0^t \int g(x)\pi_s^1(\diff x)\diff s =\int_0^t g\left(\overline{f}(s)\right)\diff s,
\]
hence, $\pi_s^1$ is in fact the Dirac measure at $\overline{f}(s)$ for $s{<}t_0^a$.

Relation~\eqref{eqZb} gives that, on the time interval $I_a$ and under the initial conditions~\eqref{InitOp}, then   the sequence of processes $(\overline{F}_{N_k}(t))$ is converging in distribution to $(\overline{f}(t))$ such that

\begin{multline}\label{f-int}
  1{-}\overline{f}(t)= 1{-}\overline{f}_0{+}\beta_6 t {-} \delta_6\rho_1\int_0^t\int_{\mathbb{R}_+}\frac{1{-}c_m{-}x}{x}\pi_s^1(\diff x)\diff s\\
  = 1{-}\overline{f}_0{+}\beta_6 t {-} \delta_6\rho_1\int_0^t\frac{1{-}c_m{-}\overline{f}(s)}{\overline{f}(s)}\diff s
\end{multline}
with Proposition~\ref{CVOcc1} and Notation~\eqref{rho}.

Hence, by uniqueness of the solution of the integral equation,
\[
\overline{f}(t)  = \overline{f}_0{-}\delta_6(\rho_6{+}\rho_1)t{+} \delta_6\rho_1(1{-}c_m)\int_0^t\frac{1}{\overline{f}(s)}\diff s
\]
holds on $I_a$, for the convergence in distribution, we thus have
\[
\lim_{N\to+\infty} \left(\overline{F}_{N}(t), t{\in}I_a\right)=\left(\overline{f}(t), t{\in}I_a\right). 
\]
The  equilibrium point of $(\overline{f}(t))$ is $\overline{f}_\infty{=}{\rho_1}(1{-}c_m)/{(\rho_6{+}\rho_1)}$,
if $\overline{f}_0{<}\overline{f}_\infty$, then $\overline{f}(t){\ge}\overline{f}_0$ for all $t{\ge}0$, and otherwise $\overline{f}(t){\ge}\overline{f}_{\infty}$. By induction, this implies that the convergence  in distribution of $(\overline{F}_{N}(t))$ can be extended on time intervals $(0,nt_0^a)$, for all $n{\ge}1$ and, consequently,  on $\R_+$. We summarize our results.
\begin{theorem}[Law of Large Numbers]\label{OLLN}
  If
  \[
  \lim_{N\to+\infty} \frac{F^0_N(0)}{N}=\overline{f}_0{\in}(0,1{-}c_m),
  \]
and $(G_N(0){=},Z_N(0){=}(m_0,z_0)$, then, for the convergence in distribution,
  \[
  \lim_{N\to+\infty} \left(\frac{F^0_N(Nt)}{N}\right)=(\overline{f}(t)),
  \]
  where $(\overline{f}(t))$ is the solution of the ODE
\begin{equation}\label{ODEf}
\frac{\diff \overline{f}}{\diff t} (t)= {-}\delta_6\left(\rho_6 {+}\rho_1\right){+}\delta_6\rho_1(1{-}c_m)\frac{1}{\overline{f}(t)},
\end{equation}
where  $\rho_6$ and $\rho_1$ are defined by Relation~\eqref{rho}.
\end{theorem}
We summarize the results obtained for the convergence in distribution of the occupation measures $(\overline{\Lambda}^0_N)$. This is an extension of Proposition~\ref{prop12}. 
\begin{corollary}\label{CoroOcc}
Under the conditions of Theorem~\ref{OLLN}, the sequence of empirical distributions $(\overline{\Lambda}^0_N)$  converge in distribution to the measure $\overline{\Lambda}^0$ such that
\[
\croc{\overline{\Lambda}^0,g}{=}\int_{\R_+}\hspace{-2mm} \E\left(g\left(s,\overline{f}(s), {\cal N}_1\left(\left[0,\rho_m\frac{c_m}{\overline{f}(s)}\right]\right), {\cal N}_2\left(\left[0,\rho_1\frac{1{-}c_m{-}\overline{f}(s)}{\overline{f}(s)}\right]\right)\right)\right)\diff s,
\]
for any  continuous function $g$ on $\R_+^2{\times}\N^2$, where ${\cal N}_1$ and ${\cal N}_2$ are two independent Poisson processes on $\R_+$ with rate $1$ and $(\overline{f}(t))$ is the solution of Relation~\eqref{ODEf} with $\overline{f}(0){=}\overline{f}_0$. 
\end{corollary}
\section{Exponential Phase}\label{ExpSec}
Throughout this section,  Conditions~\eqref{SatCond} and of exponential phase of Definition~\ref{DefPhase}  hold. 
Heuristically, if there are sufficiently many polymerases, there will be an accumulation of them in the elongation phase  of rRNAs and, therefore, the output rate of all types of rRNAs is maximal. The goal of this section is to prove precise results for this assertion. 

Under this condition, for any $1{\le}j{\le}J$, the initiation rate $\alpha_{r,j}$ of rRNA of type $j$,  is larger than $\beta_{r,j}$, the  rate at which an rRNA of type $j$ grows.

\subsection*{A coupling}
We introduce a coupling to study the occupancy of the places for transcription of rRNAs. The idea is quite simple: for $1{\le}j{\le}J$, as long as $R_j^N(t)$ is strictly less than $C_{r,j}^N$,  when $U_j^N(t){=}1$, a new polymerase is added for transcription after an exponential with parameter $\alpha_{r,j}$ and if at that time $F_N(t)$ is positive, then the variable $U_j^N(t)$ remains at $1$. See the part of transcription of rRNAs in the $Q$-Matrix of our process in Section~\ref{QMatrix}.

Otherwise, if $F_N(t){=}0$, there is a total of at least $A_N{\steq{def}}N{-}C_{r,1}^N{-}\cdots{-} C_{r,J}^N{-}J$ polymerases either in transcription of mRNAs or sequestered. If $\delta{=}\min(\eta,\beta_m)$,  the duration of time after which  there will be a free polymerase which can be accommodated by the $j$th promoter of rRNAs, with probability at least $1/J$, is stochastically bounded by an exponential random variable with parameter $2\delta A_N$.  Hence, if $F_N(t){=}0$ and $U_j^N(t){=}0$, then $U_j(t)$ returns to $1$ after a duration whose distribution is stochastically bounded by  an exponential random variable with parameter $2\delta A_N/J$.

We choose $N_0$ sufficiently large, so that
\begin{equation}\label{CoupStab}
 \frac{1}{\alpha_{r,j}}{+}\frac{J}{\delta N_0} < \frac{1}{\beta_{r,j}},\qquad \forall 1{\le}j{\le}J \text{ with }\delta{=}\min(\eta,\beta_m).
\end{equation}
We are interested in the behavior of $(Q_j^N(t)){\steq{def}}(C_{r,j}^N(t){-}R_j^N(t)), 1{\le}j{\le}J$, which measures the congestion of the transcription of the rRNAs. The above coupling shows that if $N{\ge}N_0$, it can be stochastically bounded by independent queueing processes $(\overline{Q}_j(t))$, $1{\le}j{\le}J$ characterized as follows: for $1{\le}j{\le}J$, 
\begin{itemize}
\item the arrivals of customers is a Poisson process with rate $\beta_{r,j}$.
\item The distribution of the service of a customer is the distribution of the sum of two independent exponential random variables with respective parameters $\alpha_{r,j}$ and ${\delta N_0}/{J}$.  The service will be seen as the sum of the duration of a phase $\alpha_{r,j}$ and a phase $\delta N_0/J$.  
\end{itemize}
This is an $M/G/1$ queue, see Chapter~2 of~\citet{Robert}. It has a Markovian representation as $(I_j(t),\overline{Q}_j(t))$ where $I_j(t){\in}\{1,2\}$, $I_j(t){=}1$ indicates that the customer being served is in  phase $\alpha_{r,j}$ and $I_j(t){=}2$ when it is in the phase ${\delta N_0}/{J}$. 

Under Condition~\eqref{CoupStab}, $(\overline{Q}(t)){=}((I_j(t),\overline{Q}_j(t)), 1{\le}j{\le}J)$ is a positive recurrent Mar\-kov process, since  the coordinates are independent positive recurrent Markov processes. In particular if $\overline{Q}(0){\in}(\{0,1\}{\times}\N)^J$, then
\[
\inf\left\{t{>}0: \overline{Q}(t){=}((1,0),j{=}1,\ldots,J)\right\}
\]
is almost surely finite and integrable and for any $\eps{>}0$ and $T{>}0$, there exists $K$ such that
\[
\P\left(\sup_{0\le t\le T}\max_{1\le j\le J} \overline{Q}_j(t) \ge K\right)\le \eps. 
\]
Furthermore if
\[
\overline{\tau}_j^N=\inf\{t{>}0: \overline{Q}_j(t){=}0\}, \text{ with } \overline{Q}_j(0){=}C_{r,j}^N,
\]
then it is not difficult to show, with the classical law of large numbers, that, if $i{\in}\{0,1\}$,
\[
\lim_{K\to+\infty} \frac{\E_{(i,K)}(\overline{\tau}^N_j)}{N}=
c_{r,j}\left/\left(\frac{1}{{1}/{\alpha_{r,j}}+{J}/{(\delta N_0)}}{-}\beta_{r,j}\right)\right..
\]
We have thus proved the following proposition which shows that in the exponential phase, the transcription of rRNAs is essentially congested. 
\begin{theorem}[Saturation of Transcription of rRNAs]\label{TheoSatExp}
If Conditions~\eqref{SatCond} and~\eqref{ExpCond} hold and if  $F_N(0){=}N$,  $Z_N(0){=}0$ and $(U_N(0),R_N(0)){=}(0,0)$, i.e. all polymerases are initially free,  then the variable $\tau^e_N$ defined by
\begin{equation}
\tau^e_N\steq{def} \inf\{t{>}0: R_{N,j}(t){=}C_{r,j}^N, \forall 1{\le}j{\le}J\},
\end{equation}
is almost surely finite and
\[
\sup_N\frac{\E(\tau^e_N)}{N}<{+}\infty. 
\]
For any $\eps{>}0$ and $T{>}0$, there exists $K$ such that
\begin{equation}\label{eqK}
\P_{((1,0))}\left(\sup_{0\le t\le T}\max_{1\le j\le J} C_{r,j}^N{-}R_j^N(t) \ge K \right)\le \eps
\end{equation}
\end{theorem}
The variable $\tau^e_N$  is the first time when all places for transcription of rRNAs are occupied, i.e. the first instant when this part of the system is saturated. Our proposition gives an upper bound linear in $N$ for the average value of this random variable when Condition~\eqref{ExpCond} holds.

Now we investigate the asymptotic behavior of the remaining part of the system after time $\tau^e_N$. We introduce
\[
\croc{\Lambda^F_N,g}\steq{def} \int g\left(s,F_N(s))\right)\diff s
\text{ and }
\croc{\Lambda^{0,F}_N,g}\steq{def} \int g\left(s,F^0_{A_N}(s))\right)\diff s,
\]
if $g$ is a continuous function with compact support on $\R_+{\times}\N$, where $A_N$ defined by Relation~\eqref{AN} is the number of polymerases available when transcription of rRNA is saturated. The process $(F^0_{A_N}(t))$ is the solution of Relation~\eqref{SDE1} whose initial condition is the same as the process  $(F_N(t),S_N(t),Z_N(t))$. 

\begin{lemma}[Coupling with the Auxiliary Process]\ \\
If  $(F_N(0),S_N(0),Z_N(0)){=}(f,s,z){\in}\N^3$ and if $(U_N(0),R_N(0)){=}((1,C_{r,j}^N))$, then for 
 any $g{\in}{\cal C}_c(\R_+{\times}\N)$,
  \[
  \lim_{N\to+\infty} \left|\E\left(\croc{\Lambda^{0,F}_{A_N},g}\right){-}\E\left(\croc{\Lambda^{F}_N,g}\right)\right|=0.
  \]
\end{lemma}
\begin{proof}
From Relation~\eqref{eqK}, we know that for $K$ sufficiently large, the probability of the event
  \[
  {\cal E}_K\steq{def} \left\{\sup_{t\le T} A_N-\left(N-\sum_{1}^JU_j^N(t){+}R_j^N(t)\right)\le K\right\}
  \]
  is close to $1$. 
  
Given our initial state, at time $0$ there are $A_N$ polymerases either sequestered, free or in transcription of an mRNA.
On the event ${\cal E}_K$, on the time interval $[0,T]$, there may be at most $K$ additional polymerases. Since they enter this part of the system as free, at rate at least $C_m^N{-}(N{-}C_{r,1}^N\cdots{-}C_{r,J}^N)$, they go into transcription of an mRNA. Note that, almost surely, any of these $K$ polymerases may return a finite number of times as free on $[0,T]$.  Hence, with high probability, their contribution  to the integral defining the occupation measure is arbitrarily small as $N$ gets large, and so is their impact on the random variable $(F_N(t),S_N(t),Z_N(t))$.
\end{proof}

We can now state convergence results for the number of free and sequestered polymerases. It is a direct consequence of the arguments of the proof of the last lemma and Proposition~\ref{propsub}. It shows that in this case, basically, the number of free polymerases has a Poisson distribution and the process of the number of sequestered polymerases and free {\sSs} is a positive recurrent Markov process on $\N^2$.
\begin{theorem}[Free/Sequestered Polymerases and {\sSs}]\label{TheoExpFree}
  Under Conditions~\eqref{SatCond} and~\eqref{ExpCond} and if  $(F_N(0),S_N(0),Z_N(0)){=}(f,s,z){\in}\N^3$
  and $(U_N(0),R_N(0)){=}((1,C_{r,j}^N))$, then, for the convergence in distribution,
  \[
  \lim_{N\to+\infty} \int g\left(s,F_N(s))\right)\diff s
  = \int_0^{+\infty} \E\left(g\left(u,{\cal N}_1\left(0,\rho_m\right)\right)\right)\diff u,
  \]
  for any $g{\in}{\cal C}_{c}\left(\R_+{\times}\N\right)$, where ${\cal N}_1$ is a Poisson process with rate $1$ and 
  \[
  \rho_m{=}\frac{\beta_m(1{-}c_r)}{\alpha_m(c_m{+}c_r{-}1)}, \quad c_r{\steq{def}}\sum_{j=1}^J c_{r,j}.
  \]
Furthermore,  the sequence of processes $(S_N(t),Z_N(t))$ converges in distribution for the Skorohod topology to a jump process $(S(t),Z(t))$ on $\N^2$ whose transition rates are given by
  \[
  (s,z)\longrightarrow (s,z){+}
  \begin{cases}
    (1,-1)& \lambda\rho_m z,\\
    (-1,1)& \eta s,
  \end{cases}
  \quad 
  \begin{cases}
    (0,1)& \beta_6,\\
    (0,-1)& \delta_6z.
  \end{cases}
  \]
\end{theorem}
Note that the process $(S(t),Z(t))$ is a positive recurrent Markov process. Indeed, if, for $a{>}0$,
\[
H_a(s,z)\steq{def} a s{+}z,
\]
then it is easily seen that $H_a$ is a Lyapunov function for this Markov process if $a{\in}\R_+$ is chosen so that
\[
1< a < 1{+}\frac{\delta_6}{\lambda\rho_m},
\]
see Proposition~8.14 of ~\citet{Robert}. 
\section{Stationary Phase}\label{StatSec}
Conditions~\eqref{SatCond} and  of stationary phase of Definition~\ref{DefPhase} now hold.  For any type $j$ of rRNA, the initiation rate $\alpha_{r,j}$ is less than its production rate. 

\subsection*{A coupling}
As in Section~\ref{ExpSec} we introduce a simple coupling to study the occupancy of the slots for transcription of rRNAs. Since a polymerase enters in elongation phase of an rRNA of type $j{\in}\{1,\ldots,J\}$ at rate at most $\alpha_{r,j}$, it is easy to construct a coupling with $J$ independent $M/M/1$ processes $(Q_j(t))$ with respective input rate $\alpha_{r,j}$ and service rate $\beta_{r,j}$, so that the relations
\[
R_j^N(t)\le Q_j^N(t), \quad \forall t{\ge}0, 1{\le}j{\le}J,
\]
hold. See Chapter~5 of~\citet{Robert} for example. The following proposition is a direct consequence of this coupling and the fact that, for the convergence in distribution,  the hitting time of $p$ starting from a fixed initial state is exponential with respect to $p$, for $p$ large.  See Proposition~5.16 of~\cite{Robert}
\begin{proposition}
Under Conditions~\eqref{SatCond} and~\eqref{StatCond}, and if  $F_N(0){=}N$,  $Z_N(0){=}0$ and $(U_N(0),R_N(0)){=}(0,0)$, all polymerases are initially free,  then the variable $\tau_N$ defined by
\begin{equation}
\tau^s_N\steq{def} \inf\{t{>}0: R_{N,j}(t)=0, \forall 1{\le}j{\le}J\},
\end{equation}
is almost surely finite and
\[
\sup_N\frac{\E(\tau^s_N)}{N}<{+}\infty. 
\]
\end{proposition}

\begin{lemma}
Under Condition~\eqref{StatCond} then, for any $K{>}0$,
  \[
  \lim_{N\to+\infty} \P_{(u,r)}\left(\sup_{t\le NT}\frac{R_j^N(t)}{\ln(N)^2}>K\right)=0.
  \]
\end{lemma}
\begin{proof}
This is a simple consequence of the independence of the $(Q_j(t))$ and of Proposition~5.11 of~\cite{Robert}. 
\end{proof}
The above result shows that few polymerases are in transcription of an rRNA, hence the results of Section~\ref{SecSuper} on the auxiliary process can be used, in particular Theorem~\ref{OLLN}.
\begin{theorem}[Asymptotic Behavior in Stationary Phase]\label{theostat}
Under Conditions~\eqref{SatCond} and~\eqref{StatCond}, and the initial state such that
  \[
  \lim_{N\to+\infty}\left(\frac{F_N(0)}{N},\frac{S_N(0)}{N}\right)=\left(\overline{f},\overline{s}\right) \in[0,1]^2, \text{ with } \overline{f}{+}\overline{s}{=}1-c_m,
  \]
  and $(U_N(0),R_N(0)){=}(u,r){\in}(\{0,1\}{\times}\N)^J$ then, for the convergence the sequence of processes
  \[
  \lim_{N\to+\infty} \left(\frac{F_N(t)}{N},\frac{S_N(t)}{N})\right)=\left(\overline{f}(t),1{-}c_m{-}\overline{f}(t)\right),
  \]
where $(\overline{f}(t))$ is the solution of the ODE
\begin{equation}\label{ODEf2}
\frac{\diff \overline{f}}{\diff t} (t)= {-}\delta_6\left(\rho_6 {+}\rho_1\right){+}\delta_6\rho_1(1{-}c_m)\frac{1}{\overline{f}(t)},
\end{equation}
with $\rho_1{=}{\eta}/{\lambda}$ and $\rho_6{=}{\beta_6}/{\delta_6}$.

If $g{\in}{\cal C}_c(\R_+{\times}\N)$ then, for the convergence in distribution,
\[
\lim_{N\to+\infty}\left(\int_{\R_+} g(t,Z_N(t))\diff t\right) =
 \int_{\R_+} \E\left[g\left(t, {\cal N}_1\left(\left[0,\rho_1\frac{1{-}c_m{-}\overline{f}(t)}{\overline{f}(t)}\right]\right)\right)\right]\diff t,
 \]
 where ${\cal N}_1$ is a  Poisson processes on $\R_+$ with rate $1$.
\end{theorem}
In particular, the asymptotic fraction of free  polymerases is
\[
\frac{\rho_1}{\rho_6{+}\rho_1}(1{-}c_m),
\]
and, in this state, the number of free {\sSs} has a Poisson distribution with parameter $\rho_6$.

\appendix
\section{Sub-critical Case}\label{AppSC}
It is assumed throughout this section that $c_m{>}1$ holds. We give a sketch of the proof of the averaging principle at the basis of the proof of Proposition~\ref{propsub} for the sake of completeness. The analogue of this result in the super-critical case in Section~\ref{SecSuper} is quite different and more challenging. The corresponding tightness property is less clear in this case, in particular the definition of occupation measures has to include the slow processes. The arguments of the proofs of Section~\ref{SecSuper} can be used in the same way. As it will be seen,  it is easy to show that the sequences of ``slow'' processes $({S}^0_N(t)/N)$ and $({Z}^0_N(t)/N)$ are tight. 

Recall that $\mu_N$  is the occupation measure defined by Relation~\eqref{EmpSub}. For $K{>}0$, with the same notations as in the proof of Lemma~\ref{LemMIHit}, Relation~\eqref{UpMMI} gives the inequality
\[
  \E\left(\croc{\mu_N,[0,t]{\times}[0,K]}\right)\ge \int_0^t \P(Y(Ns)\le K)\diff s =\frac{1}{N}\int_0^{Nt}\P(Y(s)\le K)\diff s.
    \]
Since $(Y(t))$ is converging in distribution to a Poisson distribution with parameter $a/b$, for any $\eps{>}0$ and $t{>}0$, there exists $K_0$ and $N_0$ such that if $K{\ge}K_0$ and $N{\ge}N_0$, then $\E\left(\croc{\mu_N,[0,t]{\times}[0,K]}\right){>}(1{-}\eps)t$. 
Lemma~1.3 and~1.4 of~\citet{Kurtz} show that  the sequence  $(\mu_N)$ of random measures is tight and any limiting point $\mu_\infty$  can be expressed as
    \[
\croc{\mu_\infty,g}= \int_{\R_+{\times}\N} g\left(u,x\right)\pi_u(\diff x)\diff u 
\]
where $(\pi_u)$ is a previsible process with values in the state space of probability distributions on $\N$.

\subsection{Proof of Proposition~\ref{propsub}}\label{proofpropsub}
By integrating Relations~\eqref{SDE2} and~\eqref{SDE3}, we obtain the identities, for $t{\ge}0$, 
\begin{align}
  {S}^0_N(t)&={S}^0_N(0){+}M_6^N(t){-}\eta \int_0^t {S}^0_N(s)\diff s{+}\lambda \int_0^t {F}^0_N(s){Z}^0_N(s)\diff s,\label{IES}\\
  {Z}^0_N(t)&={Z}^0_N(0){+}M_Z^N(t){+}\beta_6 t{-}\delta_6\int_0^t {Z}^0_N(s)\diff s\label{IEZ}\\
  & \hspace{2cm} {+}\eta \int_0^t {S}^0_N(s)\diff s{-}\lambda\int_0^t {F}^0_N(s){Z}^0_N(s)\diff s,\notag
\end{align}
where $(M_6^N(t))$ and $(M_Z^N(t))$ are martingales whose previsible increasing processes are given by
\begin{align}
\croc{M_6^N}(t)&=\eta \int_0^t {S}^0_N(s)\diff s{+}\lambda \int_0^t {F}^0_N(s){Z}^0_N(s)\diff s,\label{crocM6}\\
\croc{M_Z^N}(t)&=\beta_6 t{+}\delta_6\int_0^t {Z}^0_N(s)\diff s {+}\eta \int_0^t {S}^0_N(s)\diff s{+}\lambda\int_0^t {F}^0_N(s){Z}^0_N(s)\diff s.\label{crocZ}
\end{align}

Relations~\eqref{crocM6} and~\eqref{crocZ}, Relation~\eqref{UpMMI}, and Doob's Inequality show that, for convergence in distribution, then 
\[
\lim_{N\to+\infty} \left(\frac{M^N_6(t)}{N}\right)=\lim_{N\to+\infty} \left(\frac{M^N_Z(t)}{N}\right)=0.
\]
We note that, for $t{\ge}0$, ${S}^0_N(t){\in}[0,N]$ and $0{\le}{Z}^0_N(t){\le}N{+}{\cal P}_5((0,\beta_6){\times}(0,t])$ by Relation~\eqref{SDE3}.   Relations~\eqref{IES} and~\eqref{IEZ}, and the criterion of the modulus of continuity, see~\citet{Billingsley},  give that  the sequence of processes
$\left({{S}^0_N(t)}/{N},{{Z}^0_N(t)}/{N}\right)$ 
is tight for the convergence in distribution associated to the uniform norm on compact sets of $\R_+$.

We can therefore take a subsequence  of
$\left(\mu_N,\left({{S}^0_N(t)}/{N}\right),\left({{Z}^0_N(t)}/{N}\right)\right)$
with indices $(N_k)$ converging in distribution to $(\mu_\infty,(s(t)),(z(t)))$, where $(s(t))$ and $(z(t))$ are continuous processes.

If  $f{\in}{\cal C}_{c}\left(\N{\times}\R_+^2\right)$, Relation~\eqref{SDE1} gives the identity
    \begin{multline*}
      f\left({F}^0_{N_k}(t),\frac{{S}^0_{N_k}(t)}{{N_k}},\frac{{Z}^0_{N_k}(t)}{{N_k}}\right)=
      f\left(f_{N_k},s_{N_k},z_{N_k}\right)+M_f^{N_k}(t)\\
      +\beta_m\int_0^t\nabla_{e_1}(f)\left({X}^0_{N_k}(s)\right)\left({N_k}{-}{F}^0_{N_k}(s){-}{S}^0_{N_k}(s)\right)\diff s\\
      + \alpha_m\int_0^t\nabla_{-e_1}(f)\left({X}^0_{N_k}(s)\right)\left(C_m^{N_k}{-}{N_k}{+}{F}^0_{N_k}(s){+}{S}^0_{N_k}(s)\right){F}^0_{N_k}(s)\diff s\\
       +\lambda\int_0^t\nabla_{-e_1+\frac{e_2}{{N_k}}-\frac{e_3}{{N_k}}}(f)\left({X}^0_{N_k}(s)\right){F}^0_{N_k}(s){Z}^0_{N_k}(s)\diff s\\
       +\eta\int_0^t\nabla_{e_1-\frac{e_2}{{N_k}}+\frac{e_3}{{N_k}}}(f)\left({X}^0_{N_k}(s)\right){S}^0_{N_k}(s)\diff s\\
       +\beta_6\int_0^t\nabla_{\frac{e_3}{{N_k}}}(f)\left({X}^0_{N_k}(s)\right)\diff s
       +\delta_6\int_0^t\nabla_{-\frac{e_3}{{N_k}}}(f)\left({X}^0_{N_k}(s)\right){Z}^0_{N_k}(s)\diff s,
    \end{multline*}
with the notation $\nabla_a(f)(x){=}f(x{+}a){-}f(x)$, for $a$ and $x{\in}\N{\times}\R_+^2$. 

 With the same arguments as for the martingales $(M_S^{N_k}(t))$ and $(M_Z^{N_k}(t))$, the process $(M_f^{N_k}(t))$ is converging in distribution to $0$. 
 By dividing by ${N_k}$ the last relation, and by letting $k$ go to infinity, we get
\begin{multline*}
\int_0^t\nabla_{e_1}(f)(x,s(u),z(u))\left(\beta_m{-}\left(\beta_m-\eta\right)s(u)\right)\pi_u(\diff x) \diff u\\+\int_0^t\nabla_{-e_1}(f)(x,s(u),z(u))\left(\alpha_m(c_m{-}1{+}s(u)){+}\lambda z(u)\right)x\pi_u(\diff x)\diff u=0,
\end{multline*}
and therefore
\begin{equation}\label{eqmu}
\int_0^t \int_{\N} \Omega_{s(u),z(u)}(g)(x)\pi_u(\diff x) \diff u =0,
\end{equation}
with, for $s$, $z{\ge}0$, $s{+}z{<}1$ and $x{\in}\N$,
\begin{multline*}
\Omega_{s,z}(g)(x)=\left(\rule{0mm}{4mm}\beta_m{-}\left(\beta_m{-}\eta\right)s\right) (g(x{+}1){-}g(x))\\+\left(\rule{0mm}{4mm}\alpha_m(c_m{-}1{+}s){+}\lambda z\right)(g(x{-}1){-}g(x)),
\end{multline*}
$\Omega_{s,z}$ is the infinitesimal generator of the Markov process $(Y(t))$ of Lemma~\ref{lemMM} with $a{=}a(s,z){=}\left(\beta_m{-}\left(\beta_m{-}\eta\right)s\right)$ and $b{=}b(s,z){=}\alpha_m(c_m{-}1{+}s){+}\lambda z)$. 
From Relation~\eqref{eqmu} and with the same methods as in Section~\ref{SecSuper}, we obtain that, almost surely, 
\[
\int_0^t \int_\N g(x)\pi_u(\diff x)\diff u =\int_0^t \int_\N g(x)\pi_u(\diff x)\diff u = \int_0^t \E\left(g\left(P_{u}\right)\right)\diff u
\]
holds for all $t{>}0$ and all functions $g$ with finite support on $\N$,  where $P_{u}$ is a Poisson random variable with parameter $a(s(u),z(u))/b(s(u),z(u))$, $u{\ge}0$. 

Hence, with similar arguments as in Section~\ref{SecSuper},  for $T{\ge}0$ such that  $s(t){+}z(t){<}1$ holds  for all $t{\le}T$, we obtain that the identities
\begin{align}
&s(t)=s_0{-}\eta \int_0^t s(u)\diff u{+}\lambda  \int_0^t z(u) \frac{\beta_m{-}\left(\beta_m-\eta\right)s(u)}{\alpha_m(c_m{-}1{+}s(u)){+}\lambda z(u)}\diff u,\label{eqa1}\\
&s(t){+}z(t)=s_0{+}z_0{-}\delta_6\int_0^t z(u)\diff u.\label{eqa2},
\end{align}
hold  almost surely, for $t{\le}T$. 
From Relation~\eqref{eqa2} we obtain  that $(s(t){+}z(t))$ is  a non-increasing function, hence  $s(t){+}z(t){\le}s_0{+}z_0{<}1$,  for $t{\ge}0$, the above system has therefore a unique solution defined on $\R_+$. Since the function $(s(t){+}z(t))$ is converging at infinity, Equation~\eqref{eqa2} shows that $(s(t))$ converges at infinity too. By dividing both sides of Relations~\eqref{eqa1} and~\eqref{eqa2} by $t$ and by letting $t$ got to infinity, we deduce that both limits are zero. Proposition~\ref{propsub} is proved. 

\subsection{Proof of Proposition~\ref{propsub2}}\label{proofpropsub2}
The first assertion on the convergence of the occupation is obtained in the same way but with $(s,z){=}(0,0)$, hence for $u{\ge}0$, $s(u){=}z(u){=}0$, and the operator is
\[
\Omega_{(s(u),z(u))}(g)(x)=\beta_m (g(x{+}1){-}g(x)){+}\alpha_m(c_m{-}1)(g(x{-}1){-}g(x)).
\]
Therefore $P_{u}$ is a Poisson random variable with parameter $\rho_m$. 

Let, for $k{\ge}1$, $t_k^N$ be the $k$th jump of $(Y_N(t)){\steq{def}}(S_N(t),Z_N(t))$ when the initial state is $(s,z)$, there are four random variables  $A_i$, $i{\in}\{1,2,3,4\}$,  to trigger a change of state of $(Y_N(t))$,
\begin{enumerate}
\item $A_1^N$ is a random variable such that, for $t{\ge}0$, 
\begin{equation}\label{A1Dist}
  \P\left(A_1^N{\ge}t\mid (F_N^0(s))\right)=\exp\left({-}\lambda z\int_0^tF_N^0(s)\diff s\right);
\end{equation}
\item $A_2$, $A_3$, $A_4$ are independent exponential random variables with respective parameters $\eta s$, $\beta_6$ and $\delta_6 z$,
\end{enumerate}
and, conditionally on $(F_N^0(t))$, the random variables $A_1^N$, $A_i$, $i{\in}\{2,3,4\}$ are independent.

Relation~\eqref{A1Dist} and the convergence of the sequence $(\mu_N)$ of occupation measure  of $(F_N^0(t))$ given that $A_1^N$ is converging in distribution to an exponential distribution with parameter $\lambda z\rho_m$. 

For $t{\ge}0$,  we have 
\begin{align*}
  \P_{(s,z)}\left(\rule{0mm}{4mm}Y_N\left(t_1^N\right)=(s{+}1,z{-}1), t_1^N\ge t \right)&=\P\left(\rule{0mm}{4mm}A_1^N\ge t, A_1^N\le A_2{\wedge}A_3{\wedge}A_4\right)\\
  &=\E\left(\rule{0mm}{4mm}\ind{A_1^N{\ge}t}\exp\left({-}(\eta s{+}\beta_6{+}\delta_6 z)A_1^N\right)\right),
\end{align*}
hence,
\begin{multline*}
\lim_{N\to+\infty} \P_{(s,z)}\left(\rule{0mm}{4mm}Y_N\left(t_1^N\right)=(s{+}1,z{-}1), t_1^N{\ge}t\right)\\=\frac{\lambda\rho_m z}{(\lambda\rho_m z{+}\eta s{+}\beta_6{+}\delta_6z)} e^{{-}(\lambda\rho_m z{+}\beta_6{+}\delta_6 z{+}\eta s)t},
\end{multline*}
and this last quantity is $\P_{(s,z)}(Y(t_1)=(s{+}1,z{-}1), t_1\ge t)$, 
where $(Y(t))$  is the jump process defined in Proposition~\ref{propsub2} and  $(t_i)$ is the non-decreasing sequence of its instants of jumps. A similar convergence result is obtained in the same manner for the other possibilities for the first jump of $(Y_N(t))$. By induction, one can show that for $k{\ge}1$ and any sequence $(a_i){\in}\N^2$,
\[
  \lim_{N\to+\infty} \P\left(\rule{0mm}{4mm}Y_N(t_i^N)=a_i, 1{\le}i{\le}k\right)=\P(Y(t_i)=a_i, 1{\le}i{\le}k).
\]
We conclude the proof of the convergence by using directly the very definition of the Skorohod topology. See~\citet{Billingsley}. 
\section{Super-critical Case}\label{AppSuper}
The assumption $c_m{<}1$ holds throughout this section. Technical results used in Section~\ref{SecSuper} are presented here.

Recall that $(\overline{X}^0_N(t)){=}(\overline{F}^0_N(t),G^0_N(Nt),Z^0_N(Nt))$,
with 
\[
  G^0_N(t)\steq{def}C_m^N{-}\left(N{-}F^0_N(t){-}S^0_N(t)\right).
\]

If $f$ be a non-negative Borelian function on $\R_+{\times}\N^2$, the SDEs~\eqref{SDE1}, \eqref{SDE2}, and~\eqref{SDE3} give directly the relations
\begin{align}
f&\left(\overline{X}^0_N(t)\right)=f\left(\overline{X}^0_N(0)\right){+}M_{f,N}(t)\label{SDEf}\\
&+\lambda N\int_0^{t} \nabla_{-\frac{e_1}{N}-e_3}(f)\left(\overline{X}^0_N(s)\right) F^0_N(Ns)Z^0_N(Ns)\diff s\notag\\
&  {+}\eta N\int_0^{t}\nabla_{\frac{e_1}{N}{+}e_3}(f)\left(\overline{X}^0_N(s)\right)\left(N{-}C_m^N{+}G^0_N(Ns){-}F^0_N(Ns)\right)\diff s\notag\\
&+\alpha_m N\int_0^{t} \nabla_{-\frac{e_1}{N}-e_2}(f)\left(\overline{X}^0_N(s)\right)G^0_N(Ns)F^0_N(Ns)\diff s\notag\\
&+\beta_m N \int_0^{t} \nabla_{\frac{e_1}{N}{+}e_2}(f)\left(\overline{X}^0_N(s)\right)\left(C_m^N{-}G^0_N(Ns)\right)\diff s\notag\\
&+\beta_6 N\int_0^{t} \nabla_{e_3}(f)\left(\overline{X}^0_N(s)\right)\diff s+\delta_6 N \int_0^{t} \nabla_{-e_3}(f)\left(\overline{X}^0_N(s)\right)Z^0_N(Ns)\diff s,\notag
\end{align}
where, for $i{\in}\{1,2,3\}$, $e_i$ is the $i$th unit vector of $\R^3$,   and  $(M_{f,N}(t))$ is a local martingale and its previsible increasing process is given by
\begin{align}
  \croc{M_{f,N}}(t)&=\lambda N \int_0^{t} \nabla_{-\frac{e_1}{N}-e_3}(f)\left(\overline{X}^0_N(s)\right)^2F^0_N(Ns)Z^0_N(Ns)\diff s\label{CrocMf} \\
    &+\eta N \int_0^{t} \nabla_{\frac{e_1}{N}+e_3}(f)\left(\overline{X}^0_N(s)\right)^2 (N{-}C_m^N{+}G^0_N(Ns){-}F^0_N(Ns))\diff s\notag \\
&+\alpha_m N\int_0^{t} \nabla_{-\frac{e_1}{N}{-}{e_2}}(f)\left(\overline{X}^0_N(s)\right)^2G^0_N(Ns)F^0_N(Ns)\diff s\notag\\
&+\beta_m N \int_0^{t} \nabla_{\frac{e_1}{N}{+}{e_2}}(f)\left(\overline{X}^0_N(s)\right)^2\left(C_m^N{-}G^0_N(Ns)\right)\diff s\notag\\
&+\beta_6 N\int_0^{t}  \nabla_{e_3}(f)\left(\overline{X}^0_N(s)\right)^2\diff s\notag\\
    &+\delta_6 N\int_0^{t}\nabla_{-e_3}(f)\left(\overline{X}^0_N(s)\right)^2Z^0_N(Ns)\diff s.\notag
\end{align}

\printbibliography

\end{document}